\documentclass{amsart}
\usepackage{amssymb}
\usepackage{euscript}
\usepackage{upref}

\usepackage[all,cmtip]{xy}

\usepackage[lite,initials,shortalphabetic]{amsrefs}

\makeatletter
\newcommand{\logthispage}{\@latex@warning{Something on page \thepage
    \space needs attention}}
\makeatother

\iftrue
\newcommand{\draftdate}{Draft: \today}

\makeatletter
  \def\@oddfoot{\normalfont\scriptsize \hfill \draftdate}
  \def\@evenfoot{\normalfont\scriptsize \draftdate \hfill}
\makeatother

\fi

\newcommand{\clearemptydoublepage}%
  {\newpage
    {\pagestyle{empty}%
      \cleardoublepage
    }
}

\makeatletter
  \renewcommand{\p@enumi}{\thesubsection}
\makeatother

\newenvironment{resumeenumerate}[1]
{\begin{enumerate}
 \setcounter{enumi}{#1}
 \addtocounter{enumi}{-1}
}
{\end{enumerate}
}

\makeatletter
\newcommand{\emsection}[1]{%
  \par
  \addpenalty\@secpenalty
  \vskip 6 pt plus 9 pt
  \emph{#1.}\nobreak\enspace\ignorespaces
}
\makeatother

\newcommand{\intro}{%
  \goodbreak
  \vskip 6 pt plus 9 pt
}

\numberwithin{equation}{subsection}

\newcommand{\Period}{\rlap{\enspace .}}

\DeclareMathAlphabet{\EuRm}{U}{eur}{m}{n}
\SetMathAlphabet{\EuRm}{bold}{U}{eur}{b}{n}

\newcommand{\cat}[1]{\boldsymbol{#1}}

\newcommand{\RelCat}{\mathbf{RelCat}}
\newcommand{\RelPos}{\mathbf{RelPos}}
\newcommand{\Cat}{\mathbf{Cat}}
\newcommand{\Set}{\mathbf{Set}}

\newcommand{\Inv}{\mathit{Inv}}

\newcommand{\simp}{\mathrm{s}}

\DeclareMathOperator{\Ex}{Ex}

\DeclareMathOperator{\haut}{haut}
\DeclareMathOperator{\und}{und}
\DeclareMathOperator{\we}{we}

\newcommand{\id}{\mathrm{id}}

\DeclareMathOperator{\diag}{diag}

\newcommand{\iso}{\approx}
\newcommand{\intersect}{\cap}

\newcommand{\op}{^{\mathrm{op}}}

\newcommand{\wearrow}{\overset{\scriptscriptstyle\sim}{\to}}
\newcommand{\isoarrow}{\overset{\scriptscriptstyle\approx}{\to}}

\newcommand{\subi}{_{\textsl{i}}}
\newcommand{\subt}{_{\textsl{t}}}

\newcommand{\adj}[4]{#1\negmedspace: #2\adjarrows #3:\negmedspace #4}
\newcommand{\adjarrows}{\mathchoice{\longleftrightarrow}
  {\leftrightarrow}
  {\leftrightarrow}
  {\leftrightarrow}}

\newcommand{\pushout}[3]{#1\mathbin{\mathord{\smallcoprod}_{#2}}#3}

\newcommand{\smallcoprod}{\mathchoice{\mathbin\amalg}%
               {\mathbin\amalg}%
               {{\scriptscriptstyle\mathbin{\amalg}}}%
               {{\scriptscriptstyle\mathbin{\amalg}}}}

\newcommand{\union}{\cup}
\newcommand{\bigunion}{\bigcup}

\begin{document}

\title[Relative categories]
{Relative categories: {A}nother model for the homotopy theory of
  homotopy theories}

\author{C. Barwick}
\address{Department of Mathematics, Massachusetts Institute of
  Technology, Cambridge, MA 02139}
\email{clarkbar@gmail.com}

\author{D. M. Kan}
\address{Department of Mathematics, Massachusetts Institute of
  Technology, Cambridge, MA 02139}

\date{\today}

\maketitle

\begin{abstract} We lift Charles Rezk's complete Segal space model structure on the category of simplicial spaces to a Quillen equivalent one on the category of relative categories.
\end{abstract}

\setcounter{tocdepth}{1}
\tableofcontents

\section{Introduction}
\label{sec:Intro}

\subsection{Summary}
\label{sec:summ}

The usefulness of homotopical and (co-)homological methods in so many
parts of modern mathematics seems to be due to the following two
facts:
\begin{enumerate}
\item One often runs into what we will call \emph{relative
    categories}, i.e.\ pairs $(\cat C, \cat W)$ consisting of a
  category $\cat C$ and a subcategory $\cat W \subset \cat C$ which
  contains all the objects of $\cat C$ and of which the maps are
  called \emph{weak equivalences} because one would have liked them to
  behave like isomorphisms.
\item Such a relative category $(\cat C, \cat W)$ is in essence a
  \emph{homotopy theory} because one can not only form the
  \emph{localization} of $\cat C$ with respect to $\cat W$ (often
  called its \emph{homotopy category}) which is the category obtained
  from $\cat C$ by ``formally inverting'' all the weak equivalences,
  but one can also form the more delicate \emph{simplicial
    localization} of $\cat C$ with respect to $\cat W$, which is a
  \emph{simplicial category} (i.e,\ a category enriched over
  simplicial sets) with the same objects as $\cat C$.
\end{enumerate}

In this paper we are interested in the fact that two such relative
categories give rise to the ``same'' homotopy theory if they can be
connected by a finite zigzag of DK-\emph{equivalences}, i.e.\ weak
equivalences preserving functors which induce
\begin{itemize}
\item[---] an \emph{equivalence of categories} between their homotopy
  categories, and
\item[---] \emph{weak (homotopy) equivalences} between the simplicial sets
  involved in their simplicial localizations.
\end{itemize}
One thus can ask
\begin{enumerate}
\item whether there exists on the category $\RelCat$ of small relative categories and weak equivalence-preserving functors a model structure that is a \emph{homotopy theory of homotopy theories} in the sense that it is Quillen equivalent to the ones considered by Julie Bergner, Andr\'e Joyal, Charles Rezk, and others, and
\item whether the weak equivalences in this model structure are the DK equivalences.
\end{enumerate}
Our main result in this paper is an \emph{affirmative answer} to the first of these.

An affirmative answer to the second of these questions requires a better understanding of \emph{simplicial localization functors} \cite{BK1} and will be given in \cite{BK2}.

\subsection{Further details}
\label{sec:furdet}

Our main result consists of proving that there exists a model structure on the category
  $\RelCat$ of (small) relative categories and weak equivalence
  preserving functors between them that is Quillen equivalent to
  Charles Rezk's complete Segal structure on the category $\simp\cat
  S$ of simplicial spaces (i.e.\ bisimplicial sets) and thus \emph{is a model for the theory of $\infty$-categories} (or more precisely, $(\infty,1)$-categories).  We do this by
  showing that the Reedy model structure on $\simp\cat S$ and all its
  left Bousfield localizations (and hence in particular the just
  mentioned complete Segal structure) can be lifted to Quillen
  equivalent model structure on $\RelCat$.
  
  We also obtain for each such model structure on $\simp\cat S$ also a \emph{conjugate} model structure on $\RelCat$ with the same weak equivalences and hence the same underlying relative category as the model structure discussed above. Moreover the \emph{involution} of $\RelCat$ that sends each (small) relative category to its \emph{opposite} is a Quillen equivalence (in fact an isomorphism) between these two model structures on $\RelCat$ and models the contractible space of nontrivial auto-equivalences of theories of $(\infty,1)$-categories.
  
  The proof is basically a
  relative version of Bob Thomason's arguments that the usual model
  structure on the category of simplicial sets can be lifted to a
  Quillen equivalent model structure on the category of (small)
  categories, combined with some ideas contained in a paper he wrote together with Dana Latch and Steve Wilson.

\section{An overview}
\label{sec:Oview}

This paper consists essentially of three parts.  The first part
contains a

\subsection{Formulation of our main result}
\label{sec:MnRslt}

This will be done in the first four sections,
\S\S\ref{sec:RelCat}--\ref{sec:Main}.
\begin{enumerate}
\item In \S\ref{sec:RelCat} we introduce the category $\RelCat$ of
  (small) \emph{relative categories} and \emph{relative functors}
  between them and introduce in this category notions of
  \emph{homotopic maps} and \emph{homotopy equivalences}.

  Moreover we introduce, following Thomason, a notion of \emph{Dwyer
    maps} which are a kind of neighborhood deformation retracts with
  such cofibration-like properties (which we will verify in
  \S\ref{sec:DwyMp} (\ref{sec:RetDwMp}--\ref{sec:PrpCmpDw}) as being
  closed under retracts, pushouts, and (possibly transfinite)
  compositions.
\item In \S\ref{sec:RelPst} we consider the special case of
  \emph{relative posets} and define for them two kinds of
  \emph{subdivisions}, a \emph{terminal} one and an \emph{initial} one
  which we will denote by $\xi\subt$ and $\xi\subi$ respectively.
  Unlike what happens in the case of (ordinary) posets, these two
  subdivisions are in general not each others' opposites, but only each
  others' \emph{conjugates}.  While Thomason needed only the iteration
  of one of them we will, for reasons which will become clear in
  \S\ref{sec:DwyMp} (\ref{sec:PrpCsve}--\ref{sec:PrpMnIndDw}), need the
  \emph{composition} $\xi = \xi\subt\xi\subi$ of the two of them,
  which we will refer to as the \emph{two-fold subdivision}. Of course we could just as well have used the conjugate two-fold subdivision $\overline{\xi}=\xi\subi\xi\subt$. In that case, the \emph{opposites} of our arguments then yield a Quillen equivalent \emph{conjugate model structure} with the same weak equivalences, in which the cofibrations and fibrations are the opposites of ours.
\item \label{Int:Form} In \S\ref{sec:prelims} we develop some preliminaries needed in order to formulate our main result.
\begin{enumerate}
\item We recall what is precisely meant by \emph{lifting} a
  cofibrantly generated model structure.
\item\label{it:def} We describe the \emph{Reedy model structure} on
  the category $\simp\cat S$ of bisimplicial sets, as well as its
  \emph{left Bousfield localizations}.
\item \label{it:adj}We define two \emph{adjunctions}
  \begin{displaymath}
    \adj{K_{\xi}}{\simp\cat S}{\RelCat}{N_{\xi}}
    \qquad\text{and}\qquad
    \adj{K}{\simp\cat S}{\RelCat}{N}
  \end{displaymath}
  of which the \emph{first} is the adjunction which will allow us to
  \emph{lift} the above \eqref{it:def} model structures on $\simp\cat
  S$ to Quillen equivalent ones on $\RelCat$.
\item \label{it:Form} We also formulate a \emph{key lemma}, which
  states that the two right adjoints
  \begin{displaymath}
    N_{\xi}, N\colon \RelCat \longrightarrow \simp\cat S
  \end{displaymath}
  are naturally Reedy equivalent. At
  a crucial point (in \S\ref{sec:PrfThm}) in the proof of our main result, this key lemma enables us to use,
  instead of the functor $N_{\xi}$, the much simpler \emph{simplicial nerve} functor $N$ of Charles Rezk \cite{R} (who
  called it the \emph{classifying diagram} functor).
\end{enumerate}
\item In \S\ref{sec:Main} we state our main results and mention some of its consequences.
\begin{enumerate}
\item\label{item:ourmainresult} Our main result consists of the lifts mentioned above and hence in particular the lifts of Rezk's complete Segal model structure on $\simp\cat S$ to a Quillen equivalent one on $\RelCat$.
\item Moreover, we note that for each of the resulting model structures on $\RelCat$, there is a \emph{conjugate} model structure that is connected to it by the \emph{involution} of $\RelCat$ (\ref{sec:furdet}, I).
\item We also note that the two model structures on $\RelCat$ lifted from Rezk's complete Segal structure on $\simp\cat S$ are each models for the theory of \emph{$(\infty,1)$-categories}, and that the involution relating them models the \emph{contractible space of nontrivial auto-equivalences of the theory of $(\infty,1)$-categories}.
\item Finally, we observe, after reformulating Thomason's result in our language, that our Quillen equivalences \ref{item:ourmainresult} and Thomason's Quillen equivalences are tightly connected by a simple pair of Quillen pairs.
\end{enumerate}
\end{enumerate}

\intro
The second part of the paper consists of
\subsection{A proof of the above key lemma mentioned above}
\label{sec:IntPrf}

This will be dealt with in \S\ref{sec:HomPrsFun} and
\S\ref{sec:PrfLem}.

Thomason proved this lemma for simplicial sets by using the fact that for
every simplicial set $Y$, the natural map $Y \to \Ex Y$ \cite{K} is a
weak equivalence.  However, as we were not able to relativize this
result, we will instead relativize a quite different argument that is
contained in a paper that he wrote jointly with Dana Latch and Steve
Wilson \cite{LTW}.

In \S\ref{sec:HomPrsFun} we do the following:
\begin{enumerate}
\item \label{it:cc} We note that the category $\RelCat$ is
  \emph{closed monoidal} and that the homotopy relation in $\RelCat$
  is compatible with this closure.
\item We prove that, on \emph{finite} relative posets, the subdivision
  functor $\xi\subt$, $\xi\subi$ and $\xi$ are homotopy preserving.
\item \label{it:hom} We describe sufficient conditions on functors $\RelCat \to \simp\cat
  S$ in order that they send
  homotopic maps in $\RelCat$ to homotopic maps in $\simp\cat S$.
\end{enumerate}

Finally, in \S\ref{sec:PrfLem},
\begin{resumeenumerate}{4}
\item we use these results to relativize the
  arguments used in the paper \cite{LTW}.
\end{resumeenumerate}

\intro
The third and last part of the paper consists of
\subsection{A proof of the main result}
\label{sec:IntMn}

This will be done in \S\ref{sec:DwyMp} and \S\ref{sec:PrfThm}.

The first of these, \S\ref{sec:DwyMp}, is devoted to Dwyer maps.
\begin{enumerate}
\item \label{it:DwCl} In \ref{sec:RetDwMp}--\ref{sec:PrpCmpDw} we show
  that Dwyer maps are closed under \emph{retracts}, \emph{pushouts}
  and (possibly transfinite) \emph{compositions}.
\item \label{it:suff} In \ref{sec:PrpCsve} we describe
  \emph{sufficient conditions} on a relative inclusion of relative
  posets in order that its terminal subdivision is a Dwyer map and in
  \ref{sec:PrpBdrDw} we use this to show that if, for every pair of integers $p,q
  \ge 0$, $\Delta[p,q]$ and $\partial\Delta[p,q]$ respectively denote the
  standard $(p,q)$-bisimplex and its boundary, then the inclusion
  $\partial\Delta[p,q] \to \Delta[p,q]$ induces a relative inclusion
  \eqref{Int:Form}
  \begin{displaymath}
    K_{\xi}\partial\Delta[p,q] \longrightarrow K_{\xi}\Delta[p,q] \in
    \RelPos
  \end{displaymath}
  which is a Dwyer map.
\item In \ref{sec:PrpMnIndDw} we then use (i) and (ii) to show that every \emph{monomorphism} $L \to M \in
  \simp\cat S$ gives rise to a \emph{Dwyer map} $K_{\xi}L \to K_{\xi}M
  \in \RelCat$.
\end{enumerate}

We finally complete the proof of our main result in \S\ref{sec:PrfThm}.
\begin{resumeenumerate}{4}
\item \label{it:rltvz} In \ref{sec:KeyLem} and \ref{sec:Oview} we
  relativize another key lemma of Thomason by showing that, \emph{up
    to a weak equivalence} in the Reedy model structure on $\simp\cat
  S$
  \emph{pushouts along a Dwyer map commute with the simplicial nerve
    functor} $N$, and
\item note that, in view of the first key lemma \eqref{Int:Form}
  \emph{the same holds for the functor} $N_{\xi}$.
\item In \ref{sec:PrpUnitWE} and \ref{sec:RelPst} we deduce from this
  that \emph{the unit $1 \to N_{\xi}K_{\xi}$ of the adjunction}
  \begin{displaymath}
    \adj{K_{\xi}}{\simp\cat S}{\RelCat}{N_{\xi}}
  \end{displaymath}
  \emph{is a natural Reedy equivalence} and that \emph{a map $L \to M
    \in \simp\cat S$ is a Reedy equivalence iff the induced map
    $N_{\xi}K_{\xi}L \to N_{\xi}K_{\xi}M \in \simp\cat S$ is so}.
\item In \ref{sec:PrfofThm} we then combine these results with the
  ones of \S\ref{sec:DwyMp} to finally prove our main result.
\end{resumeenumerate}

\intro
We end with a
\subsection{Remark}
\label{sec:IntRem}

The reader may wonder why we (and Thomason) did not prove
\ref{sec:IntMn} directly, i.e.\ without using the simplicial nerve
functor $N$, as this would have avoided the need for the first key
lemma (\ref{sec:MnRslt}\ref{it:Form}).  The reason is that such a
proof would probably have been much more complicated than the present
approach, as the proof of \ref{it:rltvz} relies heavily
on the fact that the relative posets involved in the definition of the
functor $N$ all have an \emph{initial object}, something that is not
at all the case for the functor $N_{\xi}$.

\section{Relative categories}
\label{sec:RelCat}

In this section we
\begin{enumerate}
\item introduce the category $\RelCat$ of (small) \emph{relative
    categories} and \emph{relative functors} between them,
\item define a \emph{homotopy relation} on $\RelCat$, and
\item use this to describe a very useful class of relative functors
  which are a kind of \emph{neighborhood deformation retracts} and
  have \emph{cofibration-like} properties and which, following
  Thomason \cite{Th}, we will call \emph{Dwyer maps}.
\end{enumerate}

\subsection{Relative categories and functors}
\label{sec:RelCtFn}

A \textbf{relative category} will be a pair $\cat C$ consisting of
\begin{enumerate}
\item a category, called the \textbf{underlying category} and denoted
  by $\und\cat C$, and
\item a subcategory of $\cat C$, called the \textbf{category of weak
    equivalences} and denoted by $\we\cat C$, of which the maps will
  be called \textbf{weak equivalences}, which are only subject to
  the requirement that $\we\cat C$ contains all the \emph{objects} of $\cat C$
  (and hence also their identity maps).
\end{enumerate}

Similarly a \textbf{relative functor} between two relative categories
will be a \emph{weak equivalence preserving} functor and a
\textbf{relative inclusion} $\cat A \to \cat B$ will be a relative
functor such that
\begin{displaymath}
  \und\cat A \subset \und\cat B
  \qquad\text{and}\qquad
  \we\cat A = \we\cat B\intersect\cat A
\end{displaymath}

The category of the small relative categories and the relative functors
between them will be denoted by $\RelCat$. This category comes with an \emph{involution}, i.e., the automorphism
\begin{equation*}
\Inv:\RelCat\longrightarrow\RelCat
\end{equation*}
which sends each category to its opposite.

\intro Two extreme kinds of relative categories are the
\subsection{Maximal and minimal relative categories}
\label{sec:MxRlCt}

A relative category will be called
\begin{enumerate}
\item \textbf{maximal} if \emph{all} its maps are weak equivalences, and
\item \textbf{minimal} if the \emph{only} weak equivalences are the
  \emph{identity maps}.
\end{enumerate}
Given an ordinary category $\cat B$ we will denote by
\begin{displaymath}
  \hat{\cat B}
  \qquad\text{and}\qquad
  \check{\cat B}
\end{displaymath}
respectively the maximal and the minimal relative categories which
have $\cat B$ as their underlying category.

Very useful examples are, for every integer $k \ge 0$, the relative categories
\begin{displaymath}
  \hat{\cat k}
  \qquad\text{and}\qquad
  \check{\cat k}
\end{displaymath}
where $\cat k$ denotes the $k$-arrow category
\begin{displaymath}
  0 \longrightarrow \cdots \longrightarrow k
\end{displaymath}
For instance we need these to describe the 

\subsection{Homotopy relation on $\RelCat$}
\label{sec:HmRlCt}

Given two objects $\cat Y,\cat Z \in \RelCat$ and two maps $f,g\colon
\cat Y \to \cat Z \in \RelCat$, a \textbf{strict homotopy} $h\colon f
\to g$ will be a natural weak equivalence, i.e., a map
\begin{displaymath}
  h\colon \cat Y\times\hat{\cat 1} \longrightarrow \cat Z \in \RelCat
\end{displaymath}
such that \eqref{sec:MxRlCt}
\begin{displaymath}
  h(y,0) = fy
  \qquad\text{and}\qquad
  h(y,1) = gy
\end{displaymath}
for every object or map $y \in \cat Y$.  More generally, two maps
$\cat Y \to \cat Z \in \RelCat$ will be called \textbf{homotopic} if
they can be connected by a finite zigzag of such strict homotopies.

Similarly a map $f\colon \cat Y \to \cat Z \in \RelCat$ will be called
a (strict) \textbf{homotopy equivalence} if there exists a map
$f'\colon \cat Z \to \cat Y \in \RelCat$ (called a (strict)
\textbf{homotopy inverse} of $f$) such that the compositions $f'f$ and
$ff'$ are (strictly) homotopic to the identity maps of $\cat Y$ and
$\cat Z$ respectively.

A special type of such a strict homotopy equivalence is involved in
the definition of
\subsection{Strong deformation retracts}
\label{sec:StrDR}

Given a relative inclusion $\cat A \to \cat Z$ \eqref{sec:RelCtFn},
$\cat A$ will be called a \textbf{strong deformation retract} of $\cat
Z$ if there exists a \textbf{strong deformation retraction} of $\cat
Z$ onto $\cat A$, i.e.\ a pair $(r,s)$ consisting of
\begin{enumerate}
\item a map $r\colon \cat Z \to \cat A \in \RelCat$, and
\item a strict homotopy \eqref{sec:HmRlCt} $s\colon r \to 1_{\cat Z}$
  such that
\item\label{item:homotopyssdr} for every object $A \in \cat A$, $rA = A$ and $s\colon rA \to A$
  is the identity map of $A$.
\end{enumerate}

Clearly $r$ \emph{is a strict homotopy equivalence \eqref{sec:HmRlCt}
  with the inclusion $\cat A \to \cat Z$ as a strict homotopy
  inverse}.

\intro
Using these strong deformation retracts we now define an important
class of maps in $\RelCat$ called
\subsection{Dwyer maps}
\label{sec:DwyerMap}

In his construction of a model structure on the category of small
(ordinary) categories Thomason \cite{Th} introduced \emph{Dwyer maps}
which were a kind of neighborhood deformation retracts and recently
Cisinski \cite{C} noted the existence of a slightly more general and
much more convenient notion which he called \emph{pseudo-Dwyer maps}.
Our Dwyer maps will be a relative version of these pseudo-Dwyer maps
of Cisinski, i.e.:

A \textbf{Dwyer map} will be a map $\cat A' \to \cat B \in \RelCat$
which admits a (unique) factorization
\begin{displaymath}
  \cat A' \iso \cat A \longrightarrow \cat B \in \RelCat
\end{displaymath}
in which the first map is an isomorphism and the second is what we
will call a \textbf{Dwyer inclusion}, i.e.\ a relative inclusion
\eqref{sec:RelCtFn} with the following properties:
\begin{enumerate}
\item $\cat A$ is a \emph{sieve} in $\cat B$, i.e.\ if $f\colon B_{1}
  \to B_{2} \in \cat B$ and $B_{2} \in \cat A$, then $f \in \cat A$
  (or equivalently, if there exists a \emph{characteristic relative
    functor} $\alpha\colon \cat B \to \hat{\cat 1}$ such that
  $\alpha^{-1}0 = \cat A$), and
\end{enumerate}
if $Z(\cat A, \cat B)$ or just $Z\cat A$ denotes the \textbf{cosieve
  on $\cat B$ generated by $\cat A$}, i.e.\ the full relative
subcategory of $\cat B$ spanned by the objects $B \in \cat B$ for
which there exists a map $A \to B \in \cat B$ which $A \in \cat A$ (or
equivalently the smallest cosieve in $\cat B$ containing $\cat A$),
then
\begin{resumeenumerate}{2}
\item $\cat A$ is a strong deformation retract of $Z\cat A$
  \eqref{sec:StrDR}.
\end{resumeenumerate}

The usefulness of these Dwyer maps is due to the fact that, as we will
show in \S\ref{sec:DwyMp}, they have such cofibration-like properties
as being closed under retracts, pushouts and transfinite
compositions.

The definition above of a strong deformation retract, and hence also of a Dwyer map, depends on the choice of the direction of the strict homotopy $s$ in \ref{item:homotopyssdr}. The opposite choice yields the notion of a \textbf{co-Dwyer map}, i.e., a map obtained from a Dwyer map by replacing the relative categories involved by their opposites.

\section{Relative posets and their subdivisions}
\label{sec:RelPst}

An important class of relative categories consists of the relative
posets and their subdivisions (which are again relative posets).

With each relative poset $\cat P$ one can associate two subdivisions,
a ``terminal'' subdivision $\xi\subt\cat P$ and an ``initial''
subdivision $\xi\subi\cat P$.  Unlike the corresponding subdivisions
of ordinary posets, these subdivisions care in general not each others
opposites, but merely each others ``conjugates'' in the sense that
there are canonical isomorphisms
\begin{displaymath}
  (\xi\subi\cat P)\op \iso \xi\subt(\cat P\op)
  \qquad\text{or equivalently}\qquad
  (\xi\subt\cat P)\op \iso \xi\subi(\cat P\op)
\end{displaymath}
For instance, if $\cat P = \check{\cat 2}$ \eqref{sec:MxRlCt} and
$\wearrow$ indicates a weak equivalence, then
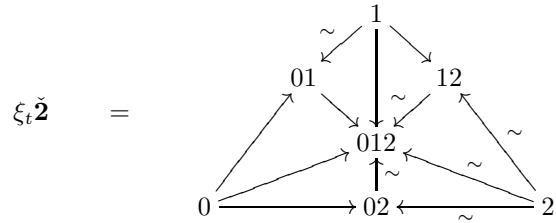
\begin{displaymath}
  \xi\subt \check{\cat 2} \qquad =\qquad
  \vcenter{\xymatrix@C=9pt@R=12pt{
      &&&{1}\ar[dr]\ar[dd]\ar[dl]_-{\sim}&&&\\
      &&{01}\ar[dr]&&{12}\ar[dl]_-{\sim}&&\\
      &&&{012}&&&\\
      {0}\ar[uurr]\ar[urrr]\ar[rrr]&&&{02}\ar[u]_-{\sim}&&&{2}\ar[uull]_-{\sim}\ar[ulll]_-{\sim}\ar[lll]^-{\sim}
    }
  }
\end{displaymath}
while
\begin{displaymath}
  \xi\subi \check{\cat 2} \qquad =\qquad
  \vcenter
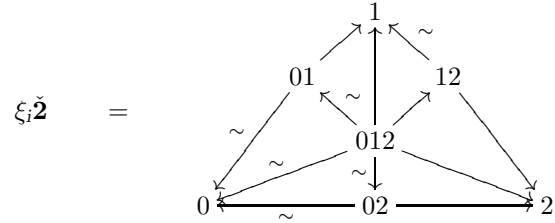
{\xymatrix@C=9pt@R=12pt{
      &&&{1}&&&\\
      &&{01}\ar[ur]\ar[ddll]_{\sim}&&{12}\ar[ul]_{\sim}\ar[ddrr]&&\\
      &&&{012}\ar[uu]\ar[ul]_-{\sim}\ar[ur]\ar[dlll]_-{\sim}\ar[d]_{\sim}\ar[drrr]&&&\\
      {0}&&&{02}\ar[lll]^-{\sim}\ar[rrr]&&&{2}
    }
  }
\end{displaymath}

\intro
In more detail:
\subsection{Relative posets}
\label{sec:RelPsts}

A \textbf{relative poset} is a relative category $\cat P$ such that
$\und\cat P$ \eqref{sec:RelCtFn} (and hence $\we\cat P$) is a poset.
The resulting full subcategory of $\RelCat$ spanned by these relative
posets will be denoted by $\RelPos$.

\subsection{Terminal and initial subdivisions}
\label{sec:TrmSbdv}

The \textbf{terminal} (resp.\ \textbf{initial}) \textbf{subdivision}
of a relative poset $\cat P$ will be the relative poset $\xi\subt\cat
P$ (resp.\ $\xi\subi\cat P$) which has
\begin{enumerate}
\item as objects the \emph{monomorphisms} \eqref{sec:MxRlCt}
  \begin{displaymath}
    \check{\cat n} \longrightarrow \cat P \in \RelPos
    \qquad(n \ge 0)
  \end{displaymath}
\item as maps
  \begin{align*}
    (x_{1}\colon \check{\cat n}_{1} \to \cat P) &\longrightarrow
    (x_{2}\colon \check{\cat n}_{2} \to \cat P)\\
    (\text{resp.}\qquad
    (x_{2}\colon \check{\cat n}_{2} \to \cat P) &\longrightarrow
    (x_{1}\colon \check{\cat n}_{1} \to \cat P))
  \end{align*}
  the commutative diagrams of the form
  \begin{displaymath}
    \xymatrix{
      {\check{\cat n}_{1}} \ar[rr] \ar[dr]_{x_{1}}
      && {\check{\cat n}_{2}} \ar[dl]^{x_{2}}\\
      & {\cat P}
    }
  \end{displaymath}
  and
\item as weak equivalences those of the above (ii) diagrams for which
  the induced map $x_{1}n_{1} \to x_{2}n_{2}$ (resp.\ $x_{2}0 \to
  x_{1}0$) is a weak equivalence in $\cat P$.
\end{enumerate}

This terminal (resp.\ initial) subdivision comes with a
\textbf{terminal} (resp.\ \textbf{initial}) \textbf{projection functor}
\begin{displaymath}
  \pi\subt\colon \xi\subt\cat P \longrightarrow \cat P
  \qquad
  (\text{resp.} \pi\subi\colon \xi\subi\cat P \longrightarrow \cat P)
\end{displaymath}
which sends an object $x\colon \check{\cat n} \to \cat P \in
\xi\subt\cat P$ (resp.\ $\xi\subi\cat P$) to the object $xn \in \cat
P$ (resp.\ $x0 \in \cat P$) and a commutative triangle as above to the
map $x_{1}n_{1} \to x_{2}n_{2} \in \cat P$ (resp.\ $x_{2}0 \to x_{1}0
\in \cat P$), which clearly implies that
\begin{resumeenumerate}{4}
\item \emph{a map in $\xi\subt\cat P$ (resp.\ $\xi\subi\cat P$) is a
    weak equivalence iff its image under $\pi\subt$ (resp.\
    $\pi\subi$) is so in $\cat P$}.
\end{resumeenumerate}

\intro
We also note the
\subsection{Naturality of the subdivisions}
\label{sec:NtSbdv}

One readily verifies that the above functions $\xi\subt$ and
$\xi\subi$ on the objects of $\RelPos$ can be extended to functors
\begin{displaymath}
  \xi\subt, \xi\subi\colon \RelPos \longrightarrow \RelPos
\end{displaymath}
by sending, for a map $f\colon \cat P \to \cat P' \in \RelPos$ every
monomorphism $\check{\cat n} \to \cat P$ to the unique monomorphism
$\check{\cat n}' \to \cat P'$ for which there exists a commutative
diagram of the form
\begin{displaymath}
  \xymatrix{
    {\check{\cat n}} \ar[r] \ar[d]
    & {\check{\cat n}'} \ar[d]\\
    {\cat P} \ar[r]_{f}
    & {\cat P'}
  }
\end{displaymath}
in  which the top map is an epimorphism.

\intro
Next we verify
\subsection{The conjugation}
\label{sec:Conj}

To verify the conjugation mentioned at the beginning of this section
we note that, using the unique isomorphisms
\begin{displaymath}
  \cat n \iso \cat n\op \qquad (n \ge 0),
\end{displaymath}
one can construct an isomorphism $\und(\xi\subi\cat P)\op \isoarrow
\und\xi\subt(\cat P\op)$ by associating with each map
\begin{displaymath}
  \vcenter{
    \xymatrix{
      {\check{\cat n}_{1}} \ar[rr] \ar[dr]_{y_{1}}
      && {\check{\cat n}_{2}} \ar[dl]^{y_{2}}\\
      & {\cat P\op}
    }
  }
  \qquad \text{in $\xi\subt(\cat P\op)$}
\end{displaymath}
the map
\begin{displaymath}
  \vcenter{
    \xymatrix{
      {\check{\cat n}_{1}} \ar[rr] \ar[d]_{\iso}
      && {\check{\cat n}_{2}} \ar[d]^{\iso}\\
      {\check{\cat n}_{1}\op} \ar[dr]_{y_{1}\op}
      && {\check{\cat n}_{2}\op} \ar[dl]^{y_{2}\op}\\
      & {\cat P}
    }
  }
  \qquad \text{in $(\xi\subi\cat P)\op$}
\end{displaymath}
A rather straightforward calculation yields that this isomorphism is
actually an isomorphism of relative posets.

\intro
We end with some
\subsection{Final comments}
\label{sec:FnlCmnt}

\begin{enumerate}
\item \label{twofold} The reason
  that, given a relative poset $\cat P$, we considered in this section
  both its terminal and its initial subdivision is that, as will be
  shown in \ref{sec:PrpCsve}--\ref{sec:PrpMnIndDw} below, in order to
  obtain the needed Dwyer maps we need the \textbf{two-fold
    subdivision} $\xi\subt\xi\subi\cat P$ and not, as one might have
  expected from Thomason's original result the iterated subdivisions
  $\xi\subt^{2}\cat P$ or $\xi\subi^{2}\cat P$.  It will therefore be
  convenient to denote the two-fold subdivision
  \begin{displaymath}
    \xi\subt \xi\subi \cat P \qquad\text{by $\xi\cat P$}
  \end{displaymath}
  and the associated composition
  \begin{displaymath}
    \xi\subt\xi\subi\cat P \xrightarrow{\enspace\pi\subt\enspace}
    \xi\subi\cat P \xrightarrow{\enspace\pi\subi\enspace} \cat P
    \qquad\text{by}\qquad\xi\cat P \xrightarrow{\enspace\pi\enspace} 
      \cat P.
  \end{displaymath}
  That Thomason did not have to do this is due to the fact that \emph{if $\cat P$ is maximal \eqref{sec:MxRlCt}, then there
    are obvious isomorphisms}
    \begin{displaymath}
      \xi\subt^{2}\cat P \iso \xi\cat P
      \qquad\text{and}\qquad
      \xi\subi^{2}\cat P \iso \xi\cat P.
    \end{displaymath}
\item\label{conjtwofold} Dually, there is a \textbf{conjugate two-fold subdivision} $\xi\subi\xi\subt$, which we denote by $\overline{\xi}\cat P$, and for which we denote the associated composition
\begin{equation*}
    \xi\subi\xi\subt\cat P \xrightarrow{\enspace\pi\subi\enspace}
    \xi\subt\cat P \xrightarrow{\enspace\pi\subt\enspace} \cat P
    \qquad\text{by}\qquad\overline{\xi}\cat P \xrightarrow{\enspace\overline{\pi}\enspace} 
      \cat P.
\end{equation*}
\item \label{convenient}
  Given a relative poset $\cat P$ it is sometimes convenient to denote
  an object
  \begin{displaymath}
    x\colon \check{\cat n} \longrightarrow \cat P \in
    \xi\subt\cat P \quad\text{or}\quad \xi\subi\cat P
  \end{displaymath}
  by the sequence
  \begin{displaymath}
    (x0, \ldots, xn)
  \end{displaymath}
  of objects of $\cat P$.
\end{enumerate}

\section{Some more preliminaries}
\label{sec:prelims}

To formulate our main result (in \ref{sec:LftQE} below) we need
\begin{enumerate}
\item a description of what is meant by \emph{lifting} a cofibrantly
  generated model structure,
\item the \emph{Reedy model structure} on the category $\simp\cat S$
  of bisimplicial sets as well as its \emph{left Bousfield
    localizations},
\item two adjunctions $\simp\cat S \adjarrows \RelCat$, and
\item a \emph{key lemma}.
\end{enumerate}

We thus start with
\subsection{Lifting model structures}
\label{sec:LiftMod}
(\cite{H}*{sec.~11.3}) Given a cofibrantly generated model
category $\cat F$ and an adjunction
\begin{displaymath}
  \adj{f}{\cat F}{\cat G}{g}
\end{displaymath}
one says that the model structure on $\cat F$ \textbf{lifts} to a
cofibrantly generated model structure on $\cat G$ if
\begin{enumerate}
\item \label{liftA} the sets of the images in $\cat G$ under the left
  adjoint $f$ of some choice of generating cofibrations and generating
  trivial cofibrations of the model structure on $\cat F$ \emph{admit
    the small object argument}, and
\item \label{liftB} the right adjoint $g$ takes all (possibly
  transfinite) compositions of pushouts of the images in $\cat G$ under
  $f$ of the generating trivial cofibrations of $\cat F$ to weak
  equivalences in $\cat F$,
\end{enumerate}
in which case
\begin{resumeenumerate}{3}
\item the \emph{generating cofibrations} and \emph{generating trivial
    cofibrations} of the model structure on $\cat G$ are the images
  under $f$ of the generating cofibrations and generating trivial
  cofibrations of the model structure on $\cat F$, and
\item \label{MpGwe} a map in $\cat G$ is a \emph{weak equivalence} or
  a \emph{fibration} iff its image under $g$ is so in $\cat F$.
\end{resumeenumerate}

\intro
Next we recall
\subsection{The Reedy model structure on $\simp\cat S$ and its left
  Bousfield localizations}
\label{sec:ReeSS}

As usual let $\cat\Delta \subset \Cat$ (the category of small categories)
be the full subcategory spanned by the posets $\cat n$ ($n \ge 0$)
\eqref{sec:MxRlCt} and let $\cat S$ and $\simp\cat S$ denote the
resulting categories
\begin{displaymath}
  \cat S = \Set^{\cat\Delta\op}
  \qquad\text{and}\qquad
  \simp\cat S = \Set^{\cat\Delta\op\times\cat\Delta\op}
\end{displaymath}
of \emph{simplicial} and \emph{bisimplicial} sets.
\begin{enumerate}
\item The \emph{standard model structure} on $\cat S$ is the
  cofibrantly generated proper model structure (\cite{H}*{pp.~210 and
    239}) in which
  \begin{enumerate}
  \item the cofibrations are the monomorphisms, and
  \item \label{WEhomeq} the weak equivalences are the maps which
    induce homotopy equivalences between the geometric realizations.
  \end{enumerate}
\item The resulting \emph{Reedy model structure} on $\simp\cat S$ is
  the cofibrantly generated proper model structure in which
  \begin{resumeenumerate}{3}
  \item \label{cofmono} the cofibrations are the monomorphisms, and
  \item\label{dfn:Reedy} the weak equivalences are the \emph{Reedy (weak)
      equivalences}, i.e.\ the maps $L \to M \in \simp\cat S$ for
    which the restrictions
    \begin{displaymath}
      L(\cat p, -) \longrightarrow M(\cat p, -) \in \cat S
      \qquad (p \ge 0)
    \end{displaymath}
    are weak equivalences (i).
  \end{resumeenumerate}
  \item \label{dfn:leftBous} A \emph{left Bousfield localization} (\cite{H}*{p.~57}) of this
    Reedy structure is any cofibrantly generated left proper model
    structure in which
    \begin{resumeenumerate}{5}
    \item the cofibrations are the monomorphisms, and
    \item the weak equivalences \emph{include} the Reedy equivalences.
    \end{resumeenumerate}
  \item \label{item:invsS} We note that the category $\simp\cat S$ admits an \emph{involution}
  \begin{equation*}
  \Inv:\simp\cat S\longrightarrow\simp\cat S,
  \end{equation*}
  which is the automorphism that sends an object $L\in\simp\cat S$ --- i.e., a functor $\cat\Delta\op\times\cat\Delta\op\longrightarrow\Set$ --- to the composition
  \begin{equation*}
  \cat\Delta\op\times\cat\Delta\op\xrightarrow{\enspace\sigma\op\times\sigma\op\enspace}\cat\Delta\op\times\cat\Delta\op\xrightarrow{\enspace L\enspace}\Set,
  \end{equation*}
  wherein $\sigma:\cat\Delta\longrightarrow\cat\Delta$ denotes the unique nontrivial automorphism of $\Delta$.
\end{enumerate}

\intro
We also need
\subsection{Two adjunctions}
\label{sec:TwoAdj}

Let $\Delta[m,n] \in \simp\cat S$ ($m,n \ge 0$) denote the
\emph{standard $(m,n)$-bisimplex} which has as its $(i,j)$-bisimplices
($i,j \ge 0$) the maps $(\cat i, \cat j) \to (\cat m, \cat n) \in
\cat\Delta\times\cat\Delta$ \eqref{sec:ReeSS}.  Our main result then
involves the adjunctions
\begin{displaymath}
  \adj{K}{\simp\cat S}{\RelCat}{N}
  \qquad\text{and}\qquad
  \adj{K_{\xi}}{\simp\cat S}{\RelCat}{N_{\xi}},
\end{displaymath}
in which $K$ and $K_{\xi}$ are the colimit preserving functors which
send $\Delta[p,q]$ ($p,q \ge 0$) to the relative categories
(\ref{sec:MxRlCt} and \ref{twofold})
\begin{displaymath}
  \check{\cat p}\times\hat{\cat q}
  \qquad\text{and}\qquad
  \xi(\check{\cat p}\times\hat{\cat q}),
\end{displaymath}
respectively, and $N$ and $N_{\xi}$ send an object $\cat X \in \RelCat$
to the bisimplicial sets which have as their $(p,q)$-bisimplices ($p,q
\ge 0$) the maps
\begin{displaymath}
  \check{\cat p}\times\hat{\cat q} \longrightarrow \cat X
  \qquad\text{and}\qquad
  \xi(\check{\cat p}\times\hat{\cat q}) \longrightarrow \cat X
  \in \RelCat,
\end{displaymath}
respectively.

The most important of these functors is the functor $N$ which Charles
Rezk called the \textbf{classifying diagram}, but which is now often
referred to as the \textbf{(simplicial) nerve functor}.  It is
connected to the functor $N_{\xi}$ by a natural transformation
\begin{displaymath}
  \pi^{*}\colon N \longrightarrow N_{\xi}
\end{displaymath}
induced by the natural transformation $\pi\colon \xi \to \id$
\eqref{twofold}.  This natural transformation $\pi^{*}$ is of
particular importance as, in view of the following key lemma
\ref{sec:keylem}, it enables us, in the proof of theorem
\ref{sec:LftQE} below, to use the functor $N$ instead of the much
more cumbersome functor $N_{\xi}$.

\subsection{A key lemma}
\label{sec:keylem}

\textbf{The natural transformation $\pi^{*}\colon N \to N_{\xi}$ is a
  natural Reedy equivalence} \eqref{sec:ReeSS}. A proof will be given in \S\S\ref{sec:HomPrsFun}-\ref{sec:PrfLem}.


\section{A statement of the main results}\label{sec:Main}

Our main result is
\subsection{Theorem. Lifting model structures on $\simp\cat S$ to Quillen equivalent ones on $\RelCat$}
\label{sec:LftQE}

\emph{The adjunction \eqref{sec:TwoAdj}}
\begin{displaymath}
  \adj{K_{\xi}}{\simp\cat S}{\RelCat}{N_{\xi}}
\end{displaymath}
\emph{lifts \eqref{sec:LiftMod} every left Bousfield localization of the Reedy model structure on
  $\simp\cat S$ (and in particular Rezk's complete Segal structure) to a Quillen equivalent cofibrantly generated left proper
  model structure on $\RelCat$ in which}
\begin{enumerate}
\item \label{LiftMi} \emph{a map is a weak equivalence iff its image
    under $N_{\xi}$ (or equivalently \eqref{sec:keylem} iff its image
    under $N$) is so in $\simp\cat S$,}
\item \label{LiftMii} \emph{a map is a fibration iff its image under
    $N_{\xi}$ is so in $\simp\cat S$,}
\item \label{LiftMiii} \emph{every cofibration is a Dwyer map
    \eqref{sec:DwyerMap},}
\item \label{LiftMiv} \emph{every cofibrant object is a relative poset
    \eqref{sec:RelPsts}.}
\end{enumerate}
\emph{Moreover, the model structure lifted from the Reedy structure itself is also right proper.}

A proof will be given in \S\ref{sec:PrfThm}.

\intro
Dualizing the proof of both \ref{sec:keylem} and \ref{sec:LftQE}, one obtains the following

\subsection{Theorem. The conjugate model structures on $\RelCat$}\label{sec:conjugatmod} \emph{The key lemma \ref{sec:keylem} and the theorem \ref{sec:LftQE} remain valid if one replaces}
\begin{enumerate}
\item \emph{$\xi$ with $\overline{\xi}$ \eqref{conjtwofold},}
\item \emph{$\pi$ with $\overline{\pi}$ \eqref{conjtwofold}, and}
\item \emph{the phrase} Dwyer map \emph{with the phrase} co-Dwyer map \eqref{sec:DwyerMap}.
\end{enumerate}

\subsection{Corollary}\label{cor:twomodstructs} \emph{The two model structures on $\RelCat$ lifted, as in \ref{sec:LftQE} and \ref{sec:conjugatmod}, from the Reedy model structure on $\simp\cat S$ or any left Bousfield localization thereof}
\begin{enumerate}
\item \emph{are Quillen equivalent,}
\item \emph{have the same weak equivalences, and hence}
\item \emph{have identical underlying relative categories.}
\end{enumerate}

\subsection{Theorem. The involution of $\RelCat$}
\begin{enumerate}
\item \emph{The involution \eqref{sec:RelCtFn}}
\begin{equation*}
\Inv:\RelCat\longrightarrow\RelCat
\end{equation*}
\emph{is an isomorphism between the two model structures of \eqref{cor:twomodstructs}.}
\item\label{item:equivsttheminvrelcat} \emph{Equivalently, a map $f\in\RelCat$ is a cofibration, fibration, or weak equivalence in one of those model structures iff $\Inv(f)\in\RelCat$ is so in the other.}
\end{enumerate}
\emsection{Proof} \ref{item:equivsttheminvrelcat} follows readily from the existence, for every pair of integers $p,q\geq 0$, of an isomorphism
\begin{equation*}
\overline{\xi}(\check{\cat p}\times\hat{\cat q})\iso(\xi(\check{\cat p}\op\times\hat{\cat q}\op))\op\iso(\xi(\check{\cat p}\times\hat{\cat q}))\op,
\end{equation*}
in which the first isomorphism is as in \ref{twofold}, and the second is induced by the \emph{unique} isomorphism
\begin{equation*}
\check{\cat p}\op\times\hat{\cat q}\op\iso\check{\cat p}\times\hat{\cat q}.
\end{equation*}

\subsection{Some $(\infty,1)$-categorical comments on the Rezk case} For the purposes of this section, let $\RelCat$ and $\simp\cat S$ denote the relative categories in which the weak equivalences are the Rezk ones, and denote by $\mathbf{RELCAT}$ the similarly defined large relative category. Then clearly
\begin{enumerate}
\item \emph{as $\simp\cat S$ is a model for the theory of $(\infty,1)$-categories, so is $\RelCat$.}
\end{enumerate}

To make a similar statement for the involution $\Inv:\RelCat\longrightarrow\RelCat$ \ref{sec:RelCtFn}, let $L^H$ denote the hammock localization of \cite{DK}. Then one can, for every relative category $\cat X$, define the \textbf{space} $\haut\cat X$ of auto-equivalences of $\cat X$ as the space which consists of the invertible components of the function space
\begin{equation*}
L^H\mathbf{RELCAT}(\cat X,\cat X).
\end{equation*}
It then follows from a result of To\"{e}n \cite[6.3]{To}, that the space $\haut\RelCat$ \emph{has two components, which are both contractible}. One of these contains the identity map of $\RelCat$, and thus the vertices of the other are the \emph{nontrivial auto-equivalences} of $\RelCat$.

Now we can state that
\begin{resumeenumerate}{2}
\item \emph{the involution $\Inv:\RelCat\longrightarrow\RelCat$ \eqref{sec:RelCtFn} is a nontrivial auto-equivalence of $\RelCat$, and hence it models the contractible space of the nontrivial auto-equivalences of theories of $(\infty,1)$-categories.}
\end{resumeenumerate}
\emsection{Proof} This follows readily from
\begin{enumerate}
\item the observation of To\"{e}n \cite[6.3]{To} that the involution $\Inv:\simp\cat S\longrightarrow\simp\cat S$ \eqref{item:invsS} is an automorphism of relative categories and is a nontrivial auto-equivalence of $\simp\cat S$, and
\item the commutativity of the diagram
\begin{equation*}
\xymatrix@C=18pt@R=18pt{
\RelCat\ar[d]_-N\ar[r]^-{\Inv}&\RelCat\ar[d]^-N\\
\simp\cat S\ar[r]_-{\Inv}&\simp\cat S.
}
\end{equation*}
\end{enumerate}

\intro
To next deal with Thomason's result \cite{Th} in our language we need
\subsection{Two more adjunctions}
\label{sec:MoreAdj}

Let $\widehat\Cat \subset \RelCat$ denote the full subcategory spanned
by the \emph{maximal} \eqref{sec:MxRlCt} relative categories.  Then
one has, corresponding to the adjunctions of \ref{sec:TwoAdj},
adjunctions
\begin{displaymath}
  \adj{k}{\cat S}{\widehat\Cat}{n}
  \qquad\text{and}\qquad
  \adj{k_{\xi}}{\cat S}{\widehat\Cat}{n_{\xi}}
\end{displaymath}
in which respectively $k$ and $k_{\xi}$ are the colimit preserving
functors which send the standard simplex $\Delta[q]$ $(q \ge 0)$ to
the maximal relative categories
\begin{displaymath}
  \hat{\cat q}
  \qquad\text{and}\qquad
  \xi\hat{\cat q}
\end{displaymath}
and $n$ and $n_{\xi}$ send an object $Y \in \widehat\Cat$ to the
simplicial sets which have as its $q$-simplices $(q \ge 0)$ the maps
\begin{displaymath}
  \hat{\cat q} \longrightarrow \cat Y
  \qquad\text{and}\qquad
  \xi\hat{\cat q} \longrightarrow \cat Y \in \widehat\Cat \subset
  \RelCat
\end{displaymath}
The functor $n\colon \widehat\Cat \to \cat S$ is the (classical)
\textbf{nerve} functor and is connected to the functor $n_{\xi}\colon
\widehat\Cat \to \cat S$ by a natural transformation
\begin{displaymath}
  \pi^{*}\colon n \longrightarrow n_{\xi}
\end{displaymath}
induced by the natural transformation $\pi\colon \xi \to \id$
\eqref{twofold}.

\intro
In our language Thomason's result then becomes
\subsection{Thomason's theorem \cite{Th}}
\label{sec:ThomThm}

\emph{The adjunction}
\begin{displaymath}
  \adj{k_{\xi}}{\cat S}{\widehat\Cat}{n_{\xi}}
\end{displaymath}
\emph{lifts \eqref{sec:LiftMod} the standard model structure on $\cat
  S$ \eqref{sec:ReeSS} to a Quillen equivalent cofibrantly generated
  proper model structure on $\widehat\Cat$ in which}
\begin{enumerate}
\item \emph{a map is a weak equivalence or a fibration iff its image
    under $n_{\xi}$ is so in $\cat S$,}
\item \emph{every cofibration is a Dwyer map \eqref{sec:DwyerMap}, and}
\item \emph{every cofibrant object is a relative poset \eqref{sec:RelPsts}.}
\end{enumerate}

Moreover
\begin{resumeenumerate}{4}
\item \emph{the natural transformation $\pi^{*}\colon n \to n_{\xi}$
    is a natural weak equivalence}
\end{resumeenumerate}
\emph{and hence}
\begin{resumeenumerate}{5}
\item \emph{a map is also a weak equivalence iff its image under the
    nerve functor $n$ is so in $\cat S$.}
\end{resumeenumerate}

\intro
We end this section pointing out a tight connection between our result
and Thomason's original one \cite{Th}.
\subsection{A tight connection between theorems \ref{sec:LftQE} and
  \ref{sec:ThomThm}}
\label{sec:Tight}

If one considers the category $\cat S$ as the subcategory of the
category $\simp\cat S$ consisting of the bisimplicial sets $L$ for
which
\begin{displaymath}
  L(\cat p,\cat q) = L(\cat 0, \cat q)
  \qquad\text{for all $p,q \ge 0$}
\end{displaymath}
then the inclusions
\begin{displaymath}
  \cat S \subset \simp\cat S
  \qquad\text{and}\qquad
  \widehat\Cat \subset \RelCat
\end{displaymath}
are the left adjoints in adjunctions
\begin{displaymath}
  \adj{i}{\cat S}{\simp\cat S}{r}
  \qquad\text{and}\qquad
  \adj{i}{\widehat\Cat}{\RelCat}{r}
\end{displaymath}
for which the units $1 \to ri$ are both the identity natural
transformations.  Then one readily verifies that
\subsection{Proposition}
\label{sec:PropSqr}

\emph{The diagram}
\begin{displaymath}
  \xymatrix@=4em{
    {\cat S} \ar@<.6ex>[r]^-{k_{\xi}} \ar@<-.6ex>[d]_{i}
    & {\widehat\Cat} \ar@<.6ex>[l]^-{n_{\xi}} \ar@<.6ex>[d]^{i}\\
    {\simp\cat S} \ar@<-.6ex>[u]_{r} \ar@<-.6ex>[r]_-{K_{\xi}}
    & {\RelCat} \ar@<-.6ex>[l]_-{N_{\xi}} \ar@<.6ex>[u]^{r}
  }
\end{displaymath}
\emph{in which the outside arrows are the left adjoints and the inside
ones the right adjoints has the following properties:}
\begin{enumerate}
\item \emph{The horizontal adjunctions are both Quillen equivalences
    (\ref{sec:LftQE} and \ref{sec:ThomThm}) and the vertical
    adjunctions are both Quillen pairs.}
\item \emph{The diagram commutes as a square of adjunctions and as a
    square of Quillen pairs.}
\end{enumerate}
\emph{Moreover}
\begin{resumeenumerate}{3}
\item \emph{$k_{\xi} = rK_{\xi}i$ and $n_{\xi} = rN_{\xi}i$.}
\end{resumeenumerate}

\section{Some homotopy preserving functors}
\label{sec:HomPrsFun}

In preparation for the proof (in \S\ref{sec:PrfLem} below) of lemma
\ref{sec:keylem} we here
\begin{enumerate}
\item note that the category $\RelCat$ is \emph{cartesian closed} and
  that the homotopy relation on $\RelCat$ is compatible with this
  cartesian closure,
\item prove that the subdivision functors (\S\ref{sec:RelPst})
  preserve homotopies between \emph{finite} relative posets and
\item describe a sufficient condition on a functor $\RelCat \to \simp\cat
  S$ \eqref{sec:ReeSS} in order that it sends homotopic maps in
  $\RelCat$ to homotopic maps in $\simp\cat S$.
\end{enumerate}

\intro
We thus start with
\subsection{Cartesian closure of $\RelCat$}
\label{sec:CCrelcat}

\emph{The category $\RelCat$ is cartesian closed}. That is \cite{M}*{Ch.~IV,
  sec.~6}, we have the following.
\begin{enumerate}
\item For every object $\cat Y \in \RelCat$, the functor
  \begin{displaymath}
    -\times\cat Y\colon \RelCat \longrightarrow \RelCat
  \end{displaymath}
  has a right adjoint $(-)^{\cat Y}$, which associates with an object
  $\cat Z \in \RelCat$ the \textbf{relative category of relative
    functors} $\cat Z^{\cat Y}$, which has
  \begin{enumerate}
  \item as \emph{objects} the maps $\cat Y \to \cat Z \in \RelCat$, and
  \item as \emph{maps} and \emph{weak equivalences} respectively the
    maps \eqref{sec:MxRlCt}
    \begin{displaymath}
      \cat Y\times\check{\cat 1} \longrightarrow \cat Z
      \qquad\text{and}\qquad
      \cat Y\times\hat{\cat 1} \longrightarrow \cat Z
      \qquad\text{$\in\RelCat$.}
    \end{displaymath}
  \end{enumerate}
\item For every three objects $\cat X$, $\cat Y$ and $\cat Z \in
  \RelCat$, there is \cite{M}*{Ch.~IV, sec.~6, Ex.~3} a natural \emph{enriched adjunction isomorphism}
  \begin{displaymath}
    \cat Z^{\cat X\times\cat Y} \iso (\cat Z^{\cat Y})^{\cat X}
    \in \RelCat,
  \end{displaymath}
  which sends
  \begin{resumeenumerate}{3}
  \item \label{CartClAdj} a map $f\colon \cat X \times\cat Y \to \cat
    Z$ to the map $g\colon \cat X \to \cat Z^{\cat Y}$, which sends an
    object $x \in \cat X$ to the map $gx\colon \cat Y \to \cat Z$
    given by $(gx)y = f(x,y)$ for every object $y \in \cat Y$, and
  \item a map
    \begin{displaymath}
      \cat X\times\cat Y\times\check{\cat 1} \longrightarrow \cat Z
      \qquad\text{(resp.\ $\cat X\times\cat Y\times\hat{\cat 1}
        \longrightarrow \cat Z$)}
    \end{displaymath}
    to the map
    \begin{displaymath}
      \cat X\times\check{\cat 1} \longrightarrow \cat Z^{\cat Y}
      \qquad\text{(resp.\ $\cat X\times\hat{\cat 1} \longrightarrow
        \cat Z^{\cat Y}$)}
    \end{displaymath}
    obtained from the obvious composition
    \begin{displaymath}
      \cat X\times\check{\cat 1}\times\cat Y \iso \cat X\times\cat
      Y\times\check{\cat 1} \longrightarrow \cat Z
    \qquad\text{(resp.\ $\cat X\times\hat{\cat 1}\times\cat Y \iso
        \cat X\times \cat Y\times \hat{\cat 1} \longrightarrow \cat
        Z$).}
    \end{displaymath}
  \end{resumeenumerate}
\end{enumerate}

\subsection{Proposition}
\label{sec:PrpHmt}

\emph{If two maps $f,g\colon \cat X \to \cat Y \in \RelCat$ are
  strictly homotopic \eqref{sec:HmRlCt}, then so are, for every object
$\cat Z \in \RelCat$ the induced maps \eqref{sec:CCrelcat}}
\begin{displaymath}
  f^{*},g^{*}\colon \cat Z^{\cat Y} \longrightarrow \cat Z^{\cat X},
\end{displaymath}
\emph{and hence, if $e\colon \cat X \to \cat Y \in \RelCat$ is a
  (strict) homotopy equivalence \eqref{sec:HmRlCt}, then so is, for
  every object $\cat Z \in \RelCat$, the induced maps}
\begin{displaymath}
  e^{*}\colon \cat Z^{\cat Y} \longrightarrow \cat Z^{\cat X}.
\end{displaymath}

\emsection{Proof} Given a strict homotopy $h\colon \cat X\times
\hat{\cat 1} \to \cat Y$, the desired strict homotopy is the map $\cat Z^{\cat Y}\times\hat{\cat 1} \to \cat Z^{\cat X}$
  which is adjoint (\ref{sec:CCrelcat}\ref{CartClAdj}) to the composition
  \begin{displaymath}
    \cat Z^{\cat Y} \xrightarrow{\enskip z^{h}\enskip} \cat Z^{(\cat
      X\times\hat{\cat 1})} \iso (\cat Z^{\cat X})^{\hat{\cat 1}}
  \end{displaymath}
  in which the isomorphism is as in \ref{sec:CCrelcat}\ref{CartClAdj}.

\subsection{Proposition}
\label{sec:PropSubd}

\emph{The subdivision functors $\xi\subt$, $\xi\subi$ and}
$\xi=\xi\subt\xi\subi$ (\S\ref{sec:RelPst})
\begin{enumerate}
\item \emph{preserve homotopies between maps from finite relative posets}
\end{enumerate}
\emph{and hence also}
\begin{resumeenumerate}{2}
\item \emph{preserve homotopy equivalences between finite relative posets.}
\end{resumeenumerate}
\emph{In particular,}
\begin{resumeenumerate}{3}
\item \label{homeq} \emph{for every pair of integers $p,q \ge 0$ all
    maps in the commutative diagram}
  \begin{displaymath}
    \xymatrix{
      {\xi(\check{\cat p}\times\hat{\cat q}) =
        \xi\subt\xi\subi(\check{\cat p}\times\hat{\cat q})}
      \ar[r]^-{\pi\subt\xi\subi} \ar[d]
      & {\xi\subi(\check{\cat p}\times\hat{\cat q})}
      \ar[r]^-{\pi\subi} \ar[d]
      & {\check{\cat p}\times\hat{\cat q}} \ar[d]\\
      {\xi\check{\cat p} = \xi\subt\xi\subi\check{\cat p}}
      \ar[r]^-{\pi\subt\xi\subi}
      & {\xi\subi\check{\cat p}} \ar[r]^-{\pi\subi}
      & {\check{\cat p}},
    }
  \end{displaymath}
  \emph{in which the vertical maps are induced by the projection
    $\check{\cat p}\times\hat{\cat q} \to \check{\cat p}$, are
    homotopy equivalences.}
\end{resumeenumerate}

\emsection{Proof}
We first deduce (iii) from (ii).

To do this we note that the map $\check{\cat p}\times\hat{\cat q} \to
\check{\cat p}$ is obviously a homotopy equivalence; hence, in view of (ii), so are the other two vertical maps.

Next we consider the commutative diagram
\begin{displaymath}
  \xymatrix{
    {\xi\subt\xi\subi\check{\cat p}} \ar[r]^-{\pi\subt\xi\subi}
    \ar[d]_{\xi\subt\pi\subi}
    & {\xi\subi\check{\cat p}} \ar[d]^{\pi\subi}\\
    {\xi\subt\check{\cat p}} \ar[r]^-{\pi\subt}
    & {\check{\cat p}},
  }
\end{displaymath}
for which one readily verifies that the maps going to $\check{\cat p}$
are homotopy equivalences with as homotopy inverses the maps which
send an object $i \in \check{\cat p}$ to the objects
\eqref{convenient}
\begin{displaymath}
  (0, \ldots, i)\in\xi\subt\check{\cat p}
  \qquad\text{and}\qquad
  (p-i, \ldots, p) \in \xi\subi\check{\cat p}
\end{displaymath}
respectively and the desired result now follows from the observation
that, in view of (ii), the map $\xi\subt\pi\subi$ is a weak equivalence
and thus so is the map $\pi\subt\xi\subi$.

Next we note that (ii) follows from (i). It thus remains to
prove (i).

To do this, it suffices to observe that, for every finite relative poset
$\cat P$, if
\begin{enumerate}
\item $n$ is the number of objects of $\cat P$ and one denotes the
  objects of $\cat P$ by the integers $1, \ldots, n$ in such a
  manner that, for every two such integers $a$ and $b$ one has $a \le
  b$, whenever there exists a map $a \to b \in \cat P$, and
\item $\cat J$ denotes the \emph{maximal} relative poset
  \eqref{sec:MxRlCt} which has $2n+1$ objects $j_{0}, \ldots, j_{2n}$
  and, for every integer $i$ with $0 \le i \le n-1$, maps
  \begin{displaymath}
    j_{2i} \longrightarrow j_{2i+1} \longleftarrow j_{2i+2},
  \end{displaymath}
\end{enumerate}
then we have the following.
  \begin{enumerate}
  \item \emph{There exists a map}
    \begin{displaymath}
      k\colon \xi\subt\cat P\times\cat J \longrightarrow
      \xi\subt(\cat P\times\hat{\cat 1}) \in \RelPos
    \end{displaymath}
    \emph{such that, in the notation of \ref{convenient}, for every
      object $(r_{1}, \ldots, r_{u}) \in \xi\subt\cat P$}
    \begin{displaymath}
      k\bigl((r_{1}, \ldots, r_{u}),j_{2n}\bigr) =
      \bigl((r_{1},0), \ldots, (r_{u},0)\bigr)
    \end{displaymath}
    \emph{and}
    \begin{displaymath}
      k\bigl((r_{1}, \ldots, r_{u}), j_{0}\bigr) =
      \bigl((r_{1},1), \ldots, (r_{u},1)\bigr).
    \end{displaymath}
    For in that case,
  \item \emph{for any two maps $f,g\colon \cat P \to \cat X \in
      \RelCat$ and strict homotopy \eqref{sec:HmRlCt}}
    \begin{displaymath}
      h\colon \cat P\times\hat{\cat 1} \longrightarrow \cat X \in
      \RelCat
    \end{displaymath}
    \emph{between them, the composition}
    \begin{displaymath}
      \xi\subt\cat P\times\cat J \xrightarrow{\enspace k\enspace}
      \xi\subt(\cat P\times\hat{\cat 1}) \xrightarrow{\enspace
        \xi\subt h\enspace} \xi\subt\cat X
    \end{displaymath}
    \emph{is a homotopy between $\xi\subt f$ and $\xi\subt q$.}
  \end{enumerate}

  A lengthy but essentially straightforward calculation (which we will
  leave to the reader) then yields that
  \begin{resumeenumerate}{3}
  \item \emph{such a map $k$ can be obtained by defining, for
      every integer $i$ with $0\leq i\leq n$ and every object $(p_{1},
    \ldots, p_{s},q_{1},\ldots,q_{t}) \in \xi\subt\cat P$ with $p_{s}<i\leq q_{1}$,}
    \begin{displaymath}
      k\bigl((p_1,\ldots,p_s,q_1,\ldots,q_t),j_{21}\bigr) = 
      \bigl((p_{1},0), \ldots, (p_{s},0), (q_{1},1),\ldots,(q_{t},1)\bigr),
    \end{displaymath}
    \emph{and, for every integer $i$ with $0 \le i \le n-p$ and object
      $(p_{1}, \ldots, p_{s}, q_{1}, \ldots, q_t) \in \xi\subt\cat
      P$ with $p_{s} < i < q_t$,}
    \begin{displaymath}
      k\bigl((p_{1}, \ldots, p_{s}, q_{1}, \ldots, q_{t}),
      j_{2i+1}\bigr) =
      \bigl((p_{1},0), \ldots, (p_{s},0), (q_{1},1), \ldots,
      (q_{t},1)\bigr)
    \end{displaymath}
    \emph{and}
    \begin{displaymath}
      k\bigl((p_{1}, \ldots, p_{s}, i, q_{1}, \ldots, q_{s}),
      j_{2i+1}\bigr) =
      \bigl((p_{1},0), \ldots, (p_{s},0), (i,0), (i,1), (q_{1},1),
      \ldots, (q_{t},1)\bigr).
    \end{displaymath}
  \end{resumeenumerate}

\intro
It thus remains to describe the needed sufficient condition on a
functor $\RelCat \to \simp\cat S$ \eqref{sec:ReeSS} in order that it
preserve homotopies, and for this we better first make clear what
exactly we will mean by
\subsection{Homotopic maps and homotopy equivalences in $\simp\cat S$}
\label{sec:HtpcS}

We will call
\begin{enumerate}
\item two maps $A \to B \in \simp\cat S$ \textbf{homotopic} if they can be
  connected by a finite zigzag of maps of the form $A\times\Delta[0,1]
  \to B \in \simp\cat S$, and
\item a map $f\colon A \to B \in \simp\cat S$ a \textbf{homotopy
    equivalence} if there exists a map $g\colon B \to A \in \simp\cat S$
  (called a \textbf{homotopy inverse} of $f$) such that the
  compositions $gf$ and $fg$ are homotopic to the identity maps of $A$
  and $B$ respectively.
\end{enumerate}

These definitions clearly imply that
\begin{resumeenumerate}{3}
\item \emph{every homotopy equivalence in $\simp\cat S$ is a Reedy
    equivalence \eqref{sec:ReeSS}.}
\end{resumeenumerate}

\intro Next, for every functor $\alpha\colon{\cat\Delta}\times{\cat\Delta}\to\RelCat$, let $N_{\alpha}\colon\RelCat\to\simp\cat S$ denote the functor that to every object $\cat X\in\RelCat$ and to every pair of integers $p,q\geq 0$ assigns the set of maps $\alpha(\cat p,\cat q)\to\cat X\in\RelCat$. Then one has:
\subsection{Proposition}
\label{sec:PrsHmt}

\emph{If $\iota\colon {\cat\Delta}\times{\cat\Delta} \to \RelCat$ \eqref{sec:ReeSS} is the functor that sends $(\cat p,\cat q)$ to $\hat{\cat q}$ ($p,q \ge 0$), and $\alpha\colon {\cat\Delta}\times{\cat\Delta} \to \RelCat$ is a functor for which there exists a natural transformation
  $\varepsilon\colon \alpha \to \iota$, then the functor
  $N_{\alpha}\colon \RelCat \to \simp\cat S$ sends homotopic maps in $\RelCat$ to homotopic maps in $\simp\cat S$ \eqref{sec:HtpcS}, and hence homotopy equivalences in $\RelCat$ to homotopy
    equivalences in $\simp\cat S$.}

\emph{This is in particular the case if}
\begin{enumerate}
\item $\alpha = \iota$ \emph{and} $\varepsilon = \id$
\end{enumerate}
\emph{and, for every pair of integers $p,q \ge 0$, if}
\begin{resumeenumerate}{2}
\item \label{AlphaProj} \emph{$\alpha(\cat p,\cat q) = \check{\cat
      p}\times\hat{\cat q}$ and $\varepsilon\cat q$ is the projection
    $\check{\cat p}\times\hat{\cat q} \to \hat{\cat q}$ ($q \ge 0$)},
  and
\item \label{alphacomp} $\alpha(\cat p,\cat q) = \xi(\check{\cat
    p}\times\hat{\cat q})$ \eqref{sec:FnlCmnt} \emph{and
    $\varepsilon\cat q$ is the composition}
  \begin{displaymath}
    \xi(\check{\cat p}\times\hat{\cat q})
    \xrightarrow{\enspace\pi\enspace}
    \check{\cat p}\times\hat{\cat q}
    \xrightarrow{\text{proj.}} \hat{\cat q}.
  \end{displaymath}
\end{resumeenumerate}

\emsection{Proof}
Given a homotopy $h\colon \cat X\times\hat{\cat 1} \to \cat Y
\in\RelCat$, the desired homotopy in $\cat S$ is the composition
\begin{displaymath}
  N_{\alpha}\cat X\times\Delta[0,1] \longrightarrow
  N_{\alpha}\cat X\times N_{\alpha}\hat{\cat 1} \iso
  N_{\alpha}(\cat X\times\hat{\cat 1})
  \xrightarrow{\enspace N_{\alpha}h\enspace}
  N_{\alpha}\cat Y
\end{displaymath}
in which the isomorphism in the middle is due to the fact that
$N_{\alpha}$ as a right adjoint preserves products, while the first
map is induced by the composition
\begin{displaymath}
  \Delta[0,1]\iso N_{\iota}\hat{\cat 1} \longrightarrow
  N_{\alpha}\hat{\cat 1},
\end{displaymath}
in which the first map is the obvious isomorphism, while the second is
induced by the natural transformation $\varepsilon\colon \alpha \to
\iota$.

\section{Proof of lemma \ref{sec:keylem}}
\label{sec:PrfLem}

To prove lemma \ref{sec:keylem} we have to show that \emph{for every object $\cat X \in \RelCat$ and integer $p
    \ge 0$, the map}
  \begin{displaymath}
    \pi_{p}^{*}\colon N\cat X(\cat p,-) \longrightarrow
    N_{\xi}\cat X(\cat p,-) \in \cat S
  \end{displaymath}
  \emph{is a weak equivalence.}

To prove this we recall that, for every pair of integers $p,q \ge 0$
\begin{align*}
  N\cat X(\cat p,\cat q) &=
  \RelCat(\check{\cat p}\times\hat{\cat q},\cat X) \qquad\text{and}\\
  N_{\xi}\cat X(\cat p,\cat q) &=
  \RelCat\bigl(\xi(\check{\cat p}\times\hat{\cat q}), \cat X\bigr)
\end{align*}
and embed the map $\pi_{p}^{*}$ in a commutative diagram
\begin{displaymath}
  \xymatrix{
    {\RelCat(\check{\cat p}\times\hat{-},\cat X)}
    \ar[r]^-{a}_-{\iso} \ar[d]_{\pi_{p}^{*}}
    & {\diag \bar{F}_{p}\cat X} \ar[r]^-{\diag f}
    & {\diag F_{p} \cat X} \ar[d]^{\diag k}\\
    {\RelCat\bigl(\xi(\check{\cat p}\times\hat{-}),\cat X\bigr)}
    \ar[r]^-{b}_-{\iso}
    & {\diag \bar G_{p}\cat X} \ar[r]^-{\diag g}
    & {\diag G_{p} \cat X}
  }
\end{displaymath}
in $\cat S$ and show that the maps $a$ and $b$ are isomorphisms and
that the other three are weak equivalences.

The bisimplicial sets $\bar F_{p}\cat X$, $F_{p}\cat X$, $G_{p}\cat
X$ and $\bar G_{p}\cat X$ and the maps between them are defined as
follows:
\begin{displaymath}
  \xymatrix{
    {\bar F_{p}\cat X(\cat q,\cat r) =
              \RelCat(\check{\cat p}\times\hat{\cat q},
      \cat X^{\hat{\cat 0}})} \ar[d]_{f}\\
    {F_{p}\cat X(\cat q,\cat r) =
             \RelCat(\check{\cat p}\times\hat{\cat q},
      \cat X^{\hat{\cat r}})}
    \ar@{}[r]|-{\overset{\text{\ref{sec:CCrelcat}\ref{CartClAdj}}}{\iso}}
    & {\RelCat(\hat{\cat r},
      \cat X^{\check{\cat p}\times\hat{\cat q}})} \ar[d]^{k}\\
    G_{p}\cat X(\cat q,\cat r) = \RelCat\bigl(\xi(\check{\cat
      p}\times\hat{\cat q}), \cat X^{\hat{\cat r}}\bigr)
    \ar@{}[r]|-{\overset{\text{\ref{sec:CCrelcat}\ref{CartClAdj}}}{\iso}}
    & {\RelCat(\hat{\cat r}, \cat X^{\xi(\check{\cat p}\times\hat{\cat
          q})})}\\
    {\bar G_{p}\cat X(\cat q,\cat r) = \RelCat\bigl(\xi(\check{\cat
        p}\times\hat{\cat q}), \cat X^{\hat{\cat 0}}\bigr),}
     \ar[u]^{g}
  }
\end{displaymath}
where $f$ and $g$ are induced by the unique maps $\hat{\cat r}\to
\hat{\cat 0}$ and $k$ is induced by the map $\pi\colon \xi(\check{\cat
p}\times\hat{\cat q}) \to \check{\cat p}\times\hat{\cat q}$
\eqref{sec:FnlCmnt}.

It then follows readily from \ref{sec:PrpHmt}, \ref{homeq} and
\ref{sec:PrsHmt} that the restrictions
\begin{displaymath}
  \text{$f(-,\cat r)$,\quad $g(-,\cat r)$,\quad and\quad $k(\cat q,-) \in \cat S$
  \qquad ($q,r \ge 0$)}
\end{displaymath}
are homotopy equivalences and hence weak equivalences.  Moreover, as
any map of bisimplicial sets that induces weak equivalences between
either their rows or their columns also induces a weak equivalence
between their diagonals, it follows that
\begin{displaymath}
  \text{$\diag f$,\quad $\diag g$,\quad and $\diag k$}
\end{displaymath}
are all weak equivalences.

Finally, to complete the proof, one notes that there are obvious
isomorphisms $a$ and $b$ which make the diagram commute.

\section{Dwyer maps}
\label{sec:DwyMp}

In preparation for the proof of theorem \ref{sec:LftQE} (in
\S\ref{sec:PrfThm} below) we here
\begin{enumerate}
\item note (in \ref{sec:RetDwMp}, \ref{sec:PrpPshDw} and
  \ref{sec:PrpCmpDw}) that Dwyer maps \eqref{sec:DwyerMap} are closed
  under retracts, pushouts and (possibly transfinite) compositions,
\item discuss (in \ref{sec:PrpCsve} and \ref{sec:PrpBdrDw}) a way
  of producing Dwyer maps which explains why our main result involves
  the \emph{two-fold} subdivision $\xi=\xi\subt\xi\subi$
  \eqref{sec:FnlCmnt} and not, as one might have expected from
  Thomason's original result \cite{Th}, the \emph{iterated} functors
  $\xi\subt^{2}$ and $\xi\subi^{2}$, and
\item use these results to show that \emph{every monomorphism}
  \begin{displaymath}
    L \longrightarrow M \in \simp\cat S \qquad\text{\eqref{sec:ReeSS}}
  \end{displaymath}
  \emph{gives rise to a Dwyer map \eqref{sec:TwoAdj}}
  \begin{displaymath}
    K_{\xi}L\longrightarrow K_{\xi}M \in \RelCat.
  \end{displaymath}
\end{enumerate}


\subsection{Proposition}
\label{sec:RetDwMp}

\emph{Every retract of a Dwyer map \eqref{sec:DwyerMap} is a Dwyer
  map.}

\emsection{Proof}
Let $\cat A \to \cat B$ be a Dwyer \emph{inclusion}
\eqref{sec:DwyerMap}, and let
\begin{displaymath}
  \xymatrix{
    {\cat A'} \ar[r] \ar[d]
    & {\cat A''} \ar[r] \ar[d]^{\bar f}
    & {\cat B'} \ar[d]^{f}\\
    {\cat A} \ar[r] \ar[d]
    & {Z\cat A} \ar[r] \ar[d]^{\bar g}
    & {\cat B} \ar[d]^{g}\\
    {\cat A'} \ar[r]
    & {\cat A''} \ar[r]
    & {\cat B'}
  }
\end{displaymath}
be a commutative diagram in which $gf = 1_{\cat B'}$, the horizontal
maps are relative inclusions \eqref{sec:RelCtFn}, and $(r,s)$ is a
strong deformation retraction \eqref{sec:StrDR} of $Z\cat A$
\eqref{sec:DwyerMap} onto $\cat A$.  Then a straightforward
calculation yields that $\cat A'$ is a sieve on $\cat B'$, that $\cat
A'' = Z\cat A'$ and that the pair $(r',s')$ where
\begin{displaymath}
  r' = \bar{g}r\bar{f}
  \qquad\text{and}\qquad
  s' = \bar{g}s\bar{f}\colon \bar{g}r\bar{f} \longrightarrow
  \bar{g}\bar{f} = 1_{\cat A''}
\end{displaymath}
is the desired strong deformation retraction of $\cat A'' = Z\cat A'$
onto $\cat A'$.

\subsection{Proposition}
\label{sec:PrpPshDw}

\emph{Let}
\begin{displaymath}
  \xymatrix{
    {\cat A} \ar[r]^-{s} \ar[d]_{i}
    & {\cat C} \ar[d]^{j}\\
    {\cat B} \ar[r]^-{t}
    & {\cat D}
  }
\end{displaymath}
\emph{be a pushout diagram in $\RelCat$ in which the map $i\colon \cat
  A \to \cat B$ is a Dwyer map \eqref{sec:DwyerMap}.  Then}
\begin{enumerate}
\item \label{DwMpPsht} \emph{the map $j\colon \cat C \to \cat D$ is a
    Dwyer map in which $Z\cat C = \pushout{Z\cat A}{\cat A}{\cat
      C}$, and}
\item \label{RlPstAll} \emph{if $\cat A$, $\cat B$ and $\cat C$ are
    relative posets, then so is $\cat D$.}
\end{enumerate}

\emph{Moreover}
\begin{resumeenumerate}{3}
\item \label{rstIso} \emph{the map $t\colon \cat B \to \cat D$
    restricts to isomorphisms}
  \begin{displaymath}
    X\cat A \iso X\cat C
    \qquad\text{and}\qquad
    X\cat A\intersect Z\cat A \iso X\cat C\intersect Z\cat C
  \end{displaymath}
  \emph{where $X\cat A \subset \cat B$ and $X\cat C \subset \cat D$
    denote the full relative subcategories spanned by the objects
    which are not in the image of $\cat A$ or $\cat C$.}
\end{resumeenumerate}

\emsection{Proof}
Assuming that the map $i\colon \cat A \to \cat B$ is a \emph{relative
  inclusion} \eqref{sec:RelCtFn} one shows that $\cat C$ is a sieve in
$\cat D$ by noting that the characteristic relative functor \eqref{sec:DwyerMap}
$\cat B \to \hat{\cat 1}$ and the map $\cat C \to \hat{\cat 1}$ which sends
all of $\cat C$ to $0$ yield a map $x\colon \cat D \to \hat{\cat 1}$ such
that $x^{-1}0 = \cat C$ and one shows in a similar manner that
$\pushout{Z\cat A}{\cat A}{\cat C}$ is a cosieve in $\cat D$.
Furthermore, the strong deformation retraction $(r,s)$ of $Z\cat A$
onto $\cat A$ induces a strong deformation retraction $(r',s')$ of
$\pushout{Z\cat A}{\cat A}{\cat C}$ onto $\cat C$ given by
\begin{alignat*}{3}
  r' = \pushout{r}{\cat A}{\cat C}&\colon&
  \pushout{Z\cat A}{\cat A}{\cat C} &\longrightarrow
  \pushout{\cat A}{\cat A}{\cat C} &&= \cat C\\
  s' = \pushout{s}{\cat A}{\cat C}&\colon&
  \pushout{r}{\cat A}{\cat C} &\longrightarrow
  \pushout{1_{Z\cat A}}{\cat A}{\cat C} &&=
  1_{\pushout{Z\cat A}{\cat A}{\cat C}}.
\end{alignat*}
This, together with the fact that $\pushout{Z\cat A}{\cat A}{\cat
  C}$ is a cosieve in $\cat D$, readily implies that
\begin{displaymath}
  \pushout{Z\cat A}{\cat A}{\cat C} = Z\cat C.
\end{displaymath}

To prove (iii) one notes that (i) the relative inclusion
\begin{displaymath}
  \hat{\cat 0} = \pushout{\cat A}{\cat A}{\hat{\cat 0}} \longrightarrow
  \pushout{\cat B}{\cat A}{\hat{\cat 0}}
\end{displaymath}
is a Dwyer map in which $Z\hat{\cat 0} = \pushout{Z\cat A}{\cat
  A}{\hat{\cat 0}}$ is obtained from $X\cat A\intersect Z\cat A$ by
adding a single object $0$ and, for every object $B \in X\cat
A\intersect Z\cat A$ a single weak equivalence $0 \to B$ and similarly
$\pushout{\cat B}{\cat A}{\hat{\cat 0}}$ is obtained from $X\cat A$ by
adding a single object $0$ and, for every object $B \in X\cat
A\intersect Z\cat A$ a single weak equivalence $0 \to B$.  Clearly
$\pushout{\cat D}{\cat C}{\hat{\cat 0}}$ admits a similar description in
terms of $X\cat C$ and $Z\cat C$ and the desired result now follows
from the observation that the map $\cat B \to \cat D$ induces an isomorphism
\begin{displaymath}
  \pushout{\cat B}{\cat A}{\hat{\cat 0}} \iso
  \pushout{\cat D}{\cat C}{\hat{\cat 0}} \Period
\end{displaymath}

Finally, to prove (ii), we note that if two objects $E,F \in \cat D$
are both in $\cat C$ or else both in $X\cat C$, then there is at most one
map between them as $\cat C$, and, in view of (iii), the relative categories $X\cat C \iso X\cat
A \subset \cat A$ are both relative posets. It thus remains to consider the
case that $E \in\cat C$ and $F \subset X\cat C$.  In that case, there
is no map $F \to E \in \cat D$ (because $\cat C$ is a sieve in $\cat
D$), and if there is a map $g\colon E \to F \in \cat D$, then $F \in
Z\cat C$ and $g = (s'F)(r'g)$; hence $g$ is unique because
$r'g\colon E \to r'F$ is in $\cat C$ and therefore unique.

\subsection{Proposition}
\label{sec:PrpCmpDw}

\emph{Every (possibly transfinite) composition of Dwyer maps is a
  Dwyer map.}

\emsection{Proof}
Assuming that all Dwyer maps involved are \emph{relative inclusions}
this follows readily from the following observations.
\begin{enumerate}
\item If $\cat A_{0} \to \cat A_{1}$ and $\cat A_{1} \to \cat A_{2}$
  are Dwyer maps with strong deformation retractions \eqref{sec:StrDR},
  \begin{center}
    \begin{tabular}{c@{\qquad}c@{\qquad}c}
      $(r_{0,1}, s_{0,1})$& and& $(r_{1,2}, s_{1,2})$\\*[1ex]
      of $Z(\cat A_{0}, \cat A_{1})$ onto $\cat A_{0}$& &
             of $Z(\cat A_{1},\cat A_{2})$ onto $\cat A_{1}$,
    \end{tabular}
  \end{center}
  then $\cat A_{0}$ is a sieve in $\cat A_{2}$, and one can obtain a
    strong deformation retraction
    \begin{displaymath}
      (r_{0,2}, s_{0,2}) \qquad\text{of $Z(\cat A_{0},\cat A_{2})$ onto
        $\cat A_{0}$}
    \end{displaymath}
    that restricts on $Z(\cat A_{0},\cat A_{1})$ to
    $(r_{0,1},s_{0,1})$ by ``composing'' the restriction $(r'_{1,2},s'_{1,2})$ of
    $(r_{1,2},s_{1,2})$ to $Z(\cat A_{0},\cat A_{1})$ with
    $(r_{0,1},s_{0,1})$, i.e.,\ by defining $(r_{0,2},s_{0,2})$ by
    \begin{displaymath}
      r_{0,2} = r_{0,1}r'_{1,2}
      \qquad\text{and}\qquad
      s_{0,2} = s'_{1,2}s_{0,1}.
    \end{displaymath}
\item If $\lambda$ is a limit ordinal, and
  \begin{displaymath}
    \cat A_{0} \longrightarrow \cdots \longrightarrow \cat A_{\beta}
    \longrightarrow \qquad (\beta \le \lambda)
  \end{displaymath}
  is a sequence of relative inclusions such that
  \begin{enumerate}
  \item for every limit ordinal $\gamma \le \lambda$,one has $\cat A_{\gamma} = \bigunion_{\alpha<\gamma}\cat A_{\alpha}$,
  \item for all $\beta<\lambda$, $\cat A_{0}$ is a sieve in $\cat
    A_{\beta}$, and
  \item there exist strong deformation retractions
    \begin{displaymath}
      (r_{0,\beta},s_{0,\beta}) \qquad
      \text{of $Z(\cat A_{0},\cat A_{\beta})$ onto $\cat A_{0}$}
    \end{displaymath}
    (one for each $\beta<\lambda$) such that, for each
    $\alpha<\beta<\lambda$, $(r_{0,\alpha},s_{0,\alpha})$ is the
    restriction of $(r_{0,\beta},s_{0,\beta})$ to $Z(\cat A_{0},\cat
    A_{\alpha})$,
  \end{enumerate}
  then \label{IsSieve} $\cat A_{0}$ is a sieve in $\cat A_{\lambda}$
    and there exists a (unique) strong deformation retraction
    $(r_{0,\lambda},s_{0,\lambda})$ of $Z(\cat A_{0},\cat
    A_{\lambda})$ onto $\cat A_{0}$ such that, for every
    $\beta<\lambda$, $(r_{0,\beta},s_{0,\beta})$ is the restriction of
    $(r_{0,\lambda},s_{0,\lambda})$ to $Z(\cat A_{0},\cat A_{\beta})$.
\end{enumerate}

\subsection{Proposition}
\label{sec:PrpCsve}

If
\begin{enumerate}
\item\label{item:catPisarelsubcat} $\cat P \to \cat Q \in \RelPos$ \emph{is a relative inclusion
    \eqref{sec:RelCtFn}, and}
\item\label{item:catPisacosieve} $\cat P$ \emph{is a cosieve in $\cat Q$ \eqref{sec:DwyerMap},}
\end{enumerate}
\emph{then the induced inclusion $\xi\subt\cat P \to \xi\subt\cat Q$
    \eqref{sec:RelPsts} is a Dwyer map \eqref{sec:DwyerMap}.}

\emsection{Proof}
For every object $(a_{0},\ldots, a_{n} \in \xi\subt\cat Q$
\eqref{convenient} either
\begin{enumerate}
\item \label{none} none of the $a_{i}$ $(0 \le i \le n)$ is in $\cat
  P$, or
\item \label{ExInt} there is (in view of \ref{item:catPisacosieve}) an integer $j$ with $0
  \le j \le n$ such that $a_{j} \in \cat P$ iff $i \ge j$, in which
  case
  \item \label{noneExIntConcl} $(a_{0}, \ldots, a_{n}) \in \xi\subt\cat
  P$ and $(a_{0}, \ldots, a_{n}) \in Z\xi\subt\cat P$.
\end{enumerate}

It now readily follows that $\xi\subt\cat P$ is a sieve in
$\xi\subt\cat Q$ and that the strong deformation retraction $(r,s)$
given by the formulas
\begin{align*}
  r(a_{0}, \ldots, a_{n}) &= (a_{j}, \ldots, a_{n}) \in \xi\subt\cat
  P\\
  s(a_{0}, \ldots, a_{n}) &= (a_{j}, \ldots, a_{n})
  \longrightarrow (a_{0}, \ldots, a_{n}) \in \xi\subt\cat Q
\end{align*}
is the desired one.

\subsection{Proposition}
\label{sec:PrpBdrDw}

\emph{For every pair of integers $p,q\ge 0$ let $\partial\Delta[p,q]
  \subset \Delta[p,q] \in \simp\cat S$ \eqref{sec:ReeSS} denote the
  largest subobject not containing its (only) non-degenerate (in both
  directions) $(p,q)$-bisimplex.  Then the inclusion $\partial\Delta[p,q]
  \to \Delta[p,q]$ induces \eqref{sec:TwoAdj} a Dwyer map}
\begin{displaymath}
  K_{\xi}\partial\Delta[p,q] \longrightarrow K_{\xi}\Delta[p,q] =
  \xi(\check{\cat p}\times\hat{\cat q}) = \xi\subt
  \xi\subi(\check{\cat p}\times\hat{\cat q}) \in \RelPos
  \qquad \eqref{sec:RelPsts}
\end{displaymath}

\emsection{Proof}
Let $K_{\xi\subi}\colon \simp\cat S \to \RelCat$ denote the colimit
preserving functor which, for every pair of integers $a,b \ge 0$,
sends $\Delta[a,b]$ to $\xi\subi(\check{\cat a}\times\hat{\cat b})$.
We show that
\begin{itemize}
\item[I] the inclusion $\partial\Delta[p,q]\to\Delta[p,q]$ induces an inclusion
\begin{displaymath}
K_{\xi_i}\partial\Delta[p,q]\to K_{\xi_i}\Delta[p,q]
\end{displaymath}
that satisfies \ref{item:catPisarelsubcat} and \ref{item:catPisacosieve}, implying that the resulting inclusion
\begin{displaymath}
\xi_tK_{\xi_i}\partial\Delta[p,q]\to \xi_tK_{\xi_i}\Delta[p,q]=K_{\xi}\Delta[p,q]
\end{displaymath}
is a Dwyer inclusion, and
\item[II] $K_{\xi}\partial\Delta[p,q]=\xi_tK_{\xi_i}\partial\Delta[p,q]$.
\end{itemize}

To show these, let $\cat D$ denote the poset that has as its objects the subcategories of $\check{\cat p}\times\hat{\cat q}$ of the form $\check{\cat a}\times\hat{\cat b}$ for which $\check{\cat a}$ and $\hat{\cat b}$ are relative subcategories of $\check{\cat p}$ and $\hat{\cat q}$, respectively, and as its morphisms the relative inclusions. One readily verifies the following.
\begin{enumerate}
\item For every pair of objects $\check{\cat a}_1\times\hat{\cat b}_1$ and $\check{\cat a}_2\times\hat{\cat b}_2\in\cat D$ for which both $\check{\cat a}_1\cap\check{\cat a}_2$ and $\hat{\cat b}_1\cap\hat{\cat b}_2$ are nonempty,
\begin{enumerate}
\item $(\check{\cat a}_1\times\hat{\cat b}_1)\cap(\check{\cat a}_2\times\hat{\cat b}_2)=(\check{\cat a}_1\cap\check{\cat a}_2)\times(\hat{\cat b}_1)\cap\hat{\cat b}_2)$
\item $\xi_i(\check{\cat a}_1\times\hat{\cat b}_1)\cap\xi_i(\check{\cat a}_2\times\hat{\cat b}_2)=\xi_i((\check{\cat a}_1\cap\check{\cat a}_2)\times(\hat{\cat b}_1)\cap\hat{\cat b}_2))$
\item $\xi(\check{\cat a}_1\times\hat{\cat b}_1)\cap\xi(\check{\cat a}_2\times\hat{\cat b}_2)=\xi((\check{\cat a}_1\cap\check{\cat a}_2)\times(\hat{\cat b}_1)\cap\hat{\cat b}_2))$.
\end{enumerate}
\item For every map $f\colon\check{\cat a}_1\times\hat{\cat b}_1\to\check{\cat a}_2\times\hat{\cat b}_2\in\cat D$,
\begin{enumerate}
\item $\xi_if$ is a relative inclusion, and $\xi_i(\check{\cat a}_1\times\hat{\cat b}_1)$ is a cosieve in $\xi(\check{\cat a}_2\times\hat{\cat b}_2)$, and
\item $\xi f$ is a relative inclusion, and $\xi(\check{\cat a}_1\times\hat{\cat b}_1)$ is a sieve in $\xi(\check{\cat a}_2\times\hat{\cat b}_2)$.
\end{enumerate}
\end{enumerate}

One verifies I above by noting, in view of (i)b and (ii)a, that $K_{\xi_i}\partial\Delta[p,q]$ is exactly the union in $\xi_{i}(\check{\cat p}\times\hat{\cat q})=K_{\xi_i}\Delta[p,q]$ of all the $\xi_i(\check{\cat a}\times\hat{\cat b})$'s, and thus the resulting inclusion $K_{\xi_i}\partial\Delta[p,q]\to K_{\xi_i}\Delta[p,q]$ satisfies \ref{item:catPisarelsubcat} and \ref{item:catPisacosieve}.

To obtain II above one first notes that, as above, the map $K_{\xi}\partial\Delta[p,q]\to K_{\xi}\Delta[p,q]$ is an inclusion, and thus the obvious map
$K_{\xi}\partial\Delta[p,q]\to\xi_tK_{\xi_i}\partial\Delta[p,q]$ is also an inclusion. It remains therefore to show that this map is onto. But this follows from the fact that, for every map $h\colon x\to y\in\xi_tK_{\xi_{i}}\partial\Delta[p,q]$, where $y$ is a monomorphism ${\cat n}\to K_{\xi_{i}}\partial\Delta[p,q]$, the object $y0\in K_{\xi_{i}}\partial\Delta[p,q]$ lies in some $\xi_i(\check{\cat a}\times\hat{\cat b})$ and hence, in view of (ii)b, the whole map $h\colon x\to y$ lies in $K_{\xi}\partial\Delta[p,q]$.

\intro
Finally we show, by combining \ref{sec:PrpBdrDw} with
\ref{sec:PrpPshDw} and \ref{sec:PrpCmpDw},
\subsection{Proposition}
\label{sec:PrpMnIndDw}

Every monomorphism $L \to M \in \simp\cat S$ induces
\eqref{sec:TwoAdj} a Dwyer map
\begin{displaymath}
  K_{\xi}L \longrightarrow K_{\xi}M \in \RelPos.
\end{displaymath}

\emsection{Proof}
Assume that $L$ is actually a subobject of $\cat M$ and denote by
$M^{n}$ ($n \ge -1$) the smallest subobject containing all
$(i,j)$-bisimplices with $i+j \le n$.  Then
\begin{enumerate}
\item $M = \bigunion_{n \ge -1} (M^{n}\union L)$ \quad and\quad
  $K_{\xi}M = \bigunion_{n \ge -1} K_{\xi}(M^{n}\union L)$.
\end{enumerate}
Furthermore if $\Delta_{n}(M,L)$ (resp.\ $\partial\Delta_{n}(M,L)$) ($n
\ge 0$) denotes the disjoint union of copies of $\Delta[i,j]$ (resp.\
$\partial\Delta[i,j]$), one for each non-degenerate (in both directions)
$(i,j)$-bisimplex with $i+j=n$ that is in $M^{n}\union L$, but not in
$M^{n-1}\union L$, then \ref{DwMpPsht} and \ref{sec:PrpBdrDw} imply:
\begin{resumeenumerate}{2}
\item The pushout diagram in $\simp\cat S$
  \begin{displaymath}
    \xymatrix{
      {\partial\Delta_{n}(M,L)} \ar[r] \ar[d]
      & {M^{n-1}\union L} \ar[d]\\
      {\Delta_{n}(M,L)} \ar[r]
      & {M^{n}\union L}
    }
  \end{displaymath}
  induces a pushout diagram in $\RelCat$
  \begin{displaymath}
    \xymatrix{
      {K_{\xi}\partial\Delta_{n}(M,L)} \ar[r] \ar[d]
      & {K_{\xi}(M^{n-1}\union L)} \ar[d]\\
      {K_{\xi}\Delta_{n}(M,L)} \ar[r]
      & {K_{\xi}(M^{n}\union L)},
    }
  \end{displaymath}
  in which the vertical maps are Dwyer maps. It therefore follows
  from (i) and \ref{sec:PrpCmpDw} that the map $K_{\xi}L \to K_{\xi}M$ is a Dwyer map as well.
\end{resumeenumerate}

That this map is in $\RelPos$, i.e.\ that $K_{\xi}M$ (and hence
$K_{\xi}L$) is a relative poset now can be shown by combining the
above for $L = \emptyset$ with \ref{RlPstAll} and the fact that every
(possibly transfinite) composition of relative inclusions of relative
posets is again a relative inclusion of relative posets.

\section{Proof of theorem \ref{sec:LftQE}}
\label{sec:PrfThm}

Before proving theorem \ref{sec:LftQE} (in \ref{sec:PrfofThm} below)
we
\begin{enumerate}
\item obtain a key lemma which states that, up to a weak equivalence
  in the Reedy model structure on $\simp\cat S$ (and hence
  in any left Bousfield localization thereof), pushing
  out along a Dwyer map commutes with applying the (simplicial) nerve
  functor $N\colon \RelCat \to \simp\cat S$ (and hence
  \eqref{sec:keylem} also the functor $N_{\xi}\colon \RelCat \to
  \simp\cat S$),
\end{enumerate}
and then
\begin{resumeenumerate}{2}
\item use this to show that the unit $1 \to N_{\xi}K_{\xi}$ of the
  adjunction
  \begin{displaymath}
    \adj{K_{\xi}}{\simp\cat S}{\RelCat}{N_{\xi}}
  \end{displaymath}
  is a natural Reedy weak equivalence, which in turn readily implies
  that a map $f\colon L \to M \in \simp\cat S$ is a weak equivalence
  in the Reedy model structure or any left Bousfield localization
  thereof iff the induced map $N_{\xi}K_{\xi}f\colon N_{\xi}K_{\xi}L
  \to N_{\xi}K_{\xi}M \in \simp\cat S$ is so.
\end{resumeenumerate}

\intro
We thus start with
\subsection{Another key lemma}
\label{sec:KeyLem}

\emph{Let}
\begin{displaymath}
  \xymatrix{
    {\cat A} \ar[r]^{s} \ar[d]_{i}
    & {\cat C} \ar[d]^{j}\\
    {\cat B} \ar[r]^{t}
    & {\cat D}
  }
\end{displaymath}
\emph{be a pushout diagram in $\RelCat$ in which the map $i\colon \cat
  A \to \cat B$ is a Dwyer map \eqref{sec:DwyerMap}.  Then, in the
  Reedy model structure on $\simp\cat S$} \emph{(and hence
  in any left Bousfield localization thereof),}
\begin{enumerate}
\item \label{IndWE} \emph{the induced map }
  \begin{displaymath}
    \pushout{N\cat B}{N\cat A}{N\cat C}\longrightarrow
    N\cat D \in \simp\cat S
  \end{displaymath}
  \emph{is a weak equivalence, and}
\item\label{item:leftpropertool} \emph{if $Ni$ is a weak equivalence, then so is $Nj$ and if
    $Ns$ is a weak equivalence, then so is $Nt$.}
\end{enumerate}

\emsection{Proof}
One readily verifies that (\ref{sec:DwyerMap} and \ref{rstIso})
\begin{displaymath}
  \text{$X\cat A$, $Z\cat A$ and $X\cat A\intersect Z\cat A$}
\end{displaymath}
are cosieves in $\cat B$ and that therefore the image of a map
$\check{\cat p}\times\hat{\cat q} \to \cat B$ ($p,q \ge 0$) is
\begin{enumerate}
\item either only in $X\cat A$,
\item or only in $Z\cat A$
\item or both in $X\cat A$ and in $Z\cat A$
\end{enumerate}
iff the image of the initial object $(0,0) \in \check{\cat
  p}\times\hat{\cat q}$ is.  It follows that $N\cat B$ and
\eqref{DwMpPsht} $N\cat D$ are pushouts
\begin{align*}
  N\cat B &\mathrel{\overset{\scriptscriptstyle b}{\iso}}
  \pushout{NX\cat A}{N(X\cat A\intersect Z\cat A)}{NZ\cat A}
  \qquad\text{and}\\
  N\cat D &\mathrel{\overset{\scriptscriptstyle d}{\iso}}
  \pushout{NX\cat C}{N(X\cat C\intersect Z\cat C)}{NZ\cat C}
\end{align*}
and that therefore the map $\pushout{N\cat B}{N\cat A}{N\cat C} \to
N\cat D$ admits a factorization
\begin{multline*}
  \pushout{N\cat B}{N\cat A}{N\cat C}
  \xrightarrow{\enspace a\enspace}
  \pushout{N\cat B}{NZ\cat A}{NZ\cat C}
  \mathrel{\overset{\scriptscriptstyle b}{\iso}}\\
  {NX\cat A}
  \mathbin{\mathord{\smallcoprod}_{N(X\cat A\intersect Z\cat A)}}
  {NZ\cat A}
  \mathbin{\mathord{\smallcoprod}_{NZ\cat A}}
  {NZ\cat C}=\\
  \pushout{NX\cat A}{N(X\cat A\intersect Z\cat A)}{NZ\cat C}
  \mathrel{\overset{\scriptscriptstyle c}{\iso}}
  \pushout{NX\cat C}{N(X\cat C\intersect Z\cat C)}{NZ\cat C}
  \mathrel{\overset{\scriptscriptstyle d}{\iso}}
  N\cat D
\end{multline*}
in which $c$ is induced by the isomorphisms of \ref{rstIso}, and $a$ is
induced by the inclusions $\cat A \to Z\cat A$ and $C \to Z\cat C$
\eqref{sec:DwyerMap}. 

Part (i) now follows from the observation that, in view of
\ref{AlphaProj} and the fact that (\ref{sec:HmRlCt}
and \ref{sec:DwyerMap}) the maps $\cat A \to Z\cat A$ and $\cat C \to
Z\cat C$ are homotopy equivalences, the induced maps
\begin{displaymath}
  N\cat A\longrightarrow NZ\cat A
  \qquad\text{and}\qquad
  N\cat C\longrightarrow NZ\cat C \in \simp\cat S
\end{displaymath}
are weak equivalences.

Furthermore the first half of (ii) is an immediate consequence of (i),
while the second half follows from (i) and the left properness of the
model structures involved.

\subsection{Corollary}
\label{sec:CorXi}

In view of lemma \ref{sec:keylem} proposition \ref{sec:KeyLem} remains
valid if one replaces everywhere the functor $N$ by $N_{\xi}$
\eqref{sec:TwoAdj}.

\subsection{Proposition}
\label{sec:PrpUnitWE}
\emph{The unit}
\begin{displaymath}
  \eta_{\xi}\colon 1 \longrightarrow N_{\xi}K_{\xi}
\end{displaymath}
\emph{of the adjunction $\adj{K_{\xi}}{\simp\cat S}{\RelCat}{N_{\xi}}$
\eqref{sec:TwoAdj} is a natural weak equivalence in the Reedy model
structure on $\simp\cat S$ (and hence also any left
Bousfield localization thereof).}

\subsection{Corollary}
\label{sec:CorWEBL}

\emph{A map $f\colon L \to M \in \simp\cat S$ is a weak equivalence in
  the Reedy model structure or any of its left Bousfield localizations
  iff the induced map $N_{\xi}K_{\xi}L \to N_{\xi}K_{\xi}M \in
  \simp\cat S$ is so.}

\emsection{Proof of \ref{sec:PrpUnitWE}}
We first show that
\begin{itemize}
\item[$(*)$] \emph{for every pair of integers $p,q \ge 0$, the map}
  \begin{displaymath}
    \eta_{\xi}\colon \Delta[p,q]\longrightarrow
    N_{\xi}K_{\xi}\Delta[p,q] \in \simp\cat S
  \end{displaymath}
  \emph{is a weak equivalence.}
\end{itemize}

This follows from the observation that, in the commutative
diagram
\begin{displaymath}
  \xymatrix{
    {\Delta[p,q]} \ar[r]^-{\eta_{\xi}} \ar[d]
    & {N_{\xi}K_{\xi}\Delta[p,q] =
      N_{\xi}\xi(\check{\cat p}\times\hat{\cat q})}
    \ar[d]^{\pi_{*}}\\
    {NK\Delta[p,q]} \ar[r]^-{\pi^{*}}
    & {N_{\xi}K\Delta[p,q] =
      N_{\xi}\xi(\check{\cat p}\times\hat{\cat q}),}
  }
\end{displaymath}
in which $\eta$ denotes the unit of the adjunction $\adj{K}{\simp\cat
  S}{\RelCat}{N}$ \eqref{sec:TwoAdj} and $\pi$ is as in \ref{twofold}.
$\eta$ is readily verified to be a Reedy equivalence, while $\pi^{*}$ and
$\pi_{*}$ are so in view of \ref{sec:keylem} and \ref{alphacomp} and \ref{homeq} respectively.

To deal with an arbitrary object $M \in \simp\cat S$ one notes that,
in the notation of the proof of \ref{sec:PrpMnIndDw},
\begin{displaymath}
  M = \bigunion_{n}M^{n}
  \qquad\text{and}\qquad
  N_{\xi}K_{\xi}M = \bigunion_{n}N_{\xi}K_{\xi}M^{n},
\end{displaymath}
and that it thus suffices to prove that
\begin{itemize}
\item[$(*)_{n}$] \emph{for every integer $n \ge 0$, the map}
  \begin{displaymath}
    \eta_{\xi}\colon M^{n}\longrightarrow N_{\xi}K_{\xi}M^{n}
    \in\simp\cat S
  \end{displaymath}
  \emph{is a weak equivalence.}
\end{itemize}

For $n = 0$ this is obvious, and we thus show that, for
$n>0$, $(*)_{n-1}$ implies $(*)_{n}$.

To do this, consider the commutative diagram, in which
$\Delta_{n}(M,\emptyset)$ and $\partial\Delta_{n}(M,\emptyset)$ are as in
the proof of \ref{sec:PrpMnIndDw},
\begin{displaymath}
  \xymatrix@C=8pt@R=12pt{
    {\partial\Delta_{n}(M,\emptyset)} \ar[rr] \ar[dd] \ar[dr]
    && {\Delta_{n}(M,\emptyset)} \ar[dr] \ar'[d] [dd]\\
    & {N_{\xi}K_{\xi}\partial\Delta_{n}(M,\emptyset)} \ar[rr] \ar[dd]
    && {N_{\xi}K_{\xi}\Delta_{n}(M,\emptyset)} \ar[dd] \ar@/^2ex/[dddr]\\
    {M^{n-1}} \ar[dr] \ar'[r] [rr]
    && {M^{n}} \ar[dr]\\
    & {N_{\xi}K_{\xi}M^{n-1}} \ar[rr] \ar@/_2ex/[rrrd]
    && {H} \ar[dr]\\
    &&&& {N_{\xi}K_{\xi}M^{n}}
  }
\end{displaymath}
in which the two squares are pushout squares and all maps are the
obvious ones.  It then follows from $(*)$ and $(*)_{n-1}$ above that
the slanted maps at the left and the top are weak equivalences and
so is therefore the map $M^{n}\to H$. The desired result now
follows from the observation that, in view of \ref{sec:PrpMnIndDw}, \ref{sec:KeyLem}, and
\ref{sec:CorXi} so is the map $H \to N_{\xi}K_{\xi}M^{n}$.

\intro
Now we are finally ready for the
\subsection{Proof of theorem \ref{sec:LftQE}}
\label{sec:PrfofThm}

\begin{enumerate}
\item \emph{The model structure.} To show that the Reedy model
  structure on $\simp\cat S$ lifts to a model structure on $\RelCat$ one has to
  verify \ref{liftA} and
  \ref{liftB}.  Clearly \ref{liftA}
  follows from the smallness of the prospective generating
  cofibrations and generating trivial cofibrations. To show that
  \ref{liftB} holds, one notes that, in view of
  \ref{sec:CorWEBL}, the right adjoint $N_{\xi}$ sends every prospective
  generating trivial cofibration to a weak equivalence in
  $\simp\cat S$ and that, in view of \ref{sec:PrpPshDw},
  \ref{sec:PrpCmpDw}, \ref{sec:PrpMnIndDw}, \ref{sec:KeyLem} and
  \ref{sec:CorXi}, the same holds for every (possibly transfinite)
  composition of pushouts of the prospective generating trivial
  cofibrations. Moreover, in view of \cite[Th. 3.3.20]{H}, all this applies also to any Bousfield localization of the Reedy structure.

  Furthermore
  \begin{enumerate}
  \item \ref{LiftMi} and \ref{LiftMii} follow from \ref{MpGwe},
  \item \ref{LiftMiii} follows from \ref{sec:RetDwMp},
    \ref{sec:PrpPshDw}, \ref{sec:PrpCmpDw} and \ref{sec:PrpMnIndDw}, and
  \item \ref{LiftMiv} follows similarly from \ref{sec:RetDwMp},
    \ref{sec:PrpPshDw}, \ref{sec:PrpCmpDw} and the fact that the colimit of
    every (possibly transfinite) sequence of monomorphisms of posets
    is again a poset.
  \end{enumerate}
\item \emph{The Quillen equivalence.} This follows readily from
  \ref{LiftMi} and \ref{sec:PrpUnitWE}.
\end{enumerate}

And finally
\begin{resumeenumerate}{3}
\item \emph{The (left) properness.} Left properness follows from
  \ref{item:leftpropertool}, and
  \ref{sec:PrpUnitWE} and the left properness of the model structures
  on $\simp\cat S$. The right properness of the model structure
  lifted from the Reedy model structure is a consequence of the right
  properness of the latter and the fact that the right adjoint
  preserves limits.
\end{resumeenumerate}



\begin{bibdiv}
  \begin{biblist}
  
    \bib{BK1}{article}{
      author={Barwick, Clark},
      author={Kan, Daniel M.},
      title={A characterization of simplicial localization functors},
      note={To appear}
    }

    \bib{BK2}{article}{
      author={Barwick, Clark},
      author={Kan, Daniel M.},
      title={In the category of relative categories the Rezk equivalences are exactly the DK equivalences},
      note={To appear}
    }

    \bib{C}{article}{
      author={Cisinski, Denis-Charles},
      title={La classe des morphismes de Dwyer n'est pas stable par
        r\'etractes}, 
      journal={Cahiers Topologie G\'eom. Diff\'erentielle Cat\'eg.},
      volume={40},
      year={1999},
      number={3},
      pages={227--231}
    }
    
    \bib{DK}{article}{
      author={Dwyer, William G.},
      author={Kan, Daniel M.},
      title={Calculating Simplicial Localizations}, 
      journal={J. Pure Appl. Algebra},
      volume={18},
      year={1980},
      number={1},
      pages={17--35}
    }

    \bib{H}{book}{
      author={Hirschhorn, Philip S.},
      title={Model Categories and Their Localizations},
      series={Math. Surveys and Monographs},
      volume={99},
      publisher={AMS},
      year={2003}
    }

    \bib{K}{article}{
      author={Kan, Daniel M.},
      title={On c. s. s. complexes},
      journal={Amer. J. Math.},
      volume={79},
      date={1957},
      pages={449--476}
    }

    \bib{LTW}{article}{
      author={Latch, Dana May},
      author={Thomason, Robert W.},
      author={Wilson, W. Stephen},
      title={Simplicial sets from categories},
      journal={Math. Z.},
      volume={164},
      year={1979},
      pages={195--214}
    }

    \bib{M}{book}{
      author={MacLane, Saunders},
      title={Categories for the working mathematician},
      note={Graduate Texts in Mathematics, Vol. 5},
      publisher={Springer-Verlag},
      place={New York},
      date={1971}
    }

    \bib{R}{article}{
      author={Rezk, Charles},
      title={A model for the homotopy theory of homotopy theory},
      journal={Trans. Amer. Math. Soc.},
      volume={353},
      year={2001},
      number={3},
      pages={973--1007}
    }

    \bib{Th}{article}{
      author={Thomason, Robert W.},
      title={$\mathbf{Cat}$ as a closed model category},
      journal={Cahiers de Topologie et G\'eometrie Differentielle},
      volume={21},
      number={3},
      year={1980},
      pages={305--324}
    }
    
    \bib{To}{article}{
      author={To{\"e}n, Bertrand},
      title={Vers une axiomatisation de la th\'eorie des cat\'egories
              sup\'erieures},
      journal={$K$-Theory},
      volume={34},
      year={2005},
      number={3},
      pages={233--263}
    }

  \end{biblist}
\end{bibdiv}
\end{document}